\documentclass[reqno]{amsart}

\usepackage{amssymb,epic,eepic}

\usepackage[cmtip,matrix,arrow]{xy}
\CompileMatrices
\newdir{((}{{\lhook\kern-.1em}}
\newdir{))}{{\rhook\kern-.1em}}
\newdir{ >}{{}*!/-5pt/@{>}}

\numberwithin{equation}{section}
\newcommand\note[1]%
{$^\dagger$\marginpar{\footnotesize{$^\dagger${#1}}}} 

\newtheorem{theorem}{Theorem}[section]
\newtheorem{lemma}[theorem]{Lemma}
\newtheorem{proposition}[theorem]{Proposition}
\newtheorem{corollary}[theorem]{Corollary}

\newtheorem{addendum}[theorem]{Addendum}

\theoremstyle{definition}
\newtheorem{definition}[theorem]{Definition}

\theoremstyle{remark}
\newtheorem{remark}[theorem]{Remark}
\newtheorem{example}[theorem]{Example}

\newcommand\eu{\mathfrak} 
\newcommand\lie{\mathfrak} 
\newcommand\bb{\mathbb} 
\newcommand\N{\bb{N}}
\newcommand\Z{\bb{Z}} 
\newcommand\Q{\bb{Q}}
\newcommand\R{\bb{R}} 
\newcommand\C{\bb{C}}
\renewcommand\P{{\bb P}}
\newcommand\ca{\mathcal} 
\newcommand\hw{{\ca C}}
\renewcommand\O{{\ca O}}
\newcommand\Lie{{\ca L}}
\newcommand\eps{\varepsilon}

\renewcommand\Im{\operatorname{Im}}
\newcommand\grad{\operatorname{grad}}
\newcommand\ad{\operatorname{ad}}

\newcommand\id{\operatorname{id}}
\newcommand\Hom{\operatorname{Hom}} 
\newcommand\spec{\operatorname{Spec}} 
\newcommand\proj{\operatorname{Proj}} 
\newcommand\cone{\operatorname{cone}}
\newcommand\hull{\operatorname{hull}} 

\newcommand\SU{\operatorname{SU}}

\newcommand\qu{/\kern-.7ex/} 
\newcommand\bigqu{\big/\kern-.85ex\big/} 

\newcommand\co{^{\C}} 
\newcommand\sst{^{\mathrm s\mathrm s}} 
\newcommand\sq{\sqrt{-1}}
\newcommand\inv{^{-1}} 
\newcommand\const{\frac\sq{2\pi}} 
\newcommand\tplus{\lie t^*_+} 
\newcommand\lplus{\Lambda^*_+} 
\newcommand\red{^{\text{\rm red}}} 

\title{Convexity Properties of the Moment Mapping Re-examined} 
\author{Reyer Sjamaar}
\address{School of Mathematics, Institute for Advanced Study,
Princeton, New Jersey 08540}
\curraddr{Department of Mathematics, Cornell University, Ithaca, New
York 14853-7901}
\email{sjamaar@math.cornell.edu}
\date{December 1994.  Revised July 1997.  To appear in Adv.\ in Math}

\begin{document}


\begin{abstract}
Consider a Hamiltonian action of a compact Lie group on a compact
symplectic manifold.  A theorem of Kirwan's says that the image of the
momentum mapping intersects the positive Weyl chamber in a convex
polytope.  I present a new proof of Kirwan's theorem, which gives
explicit information on how the vertices of the polytope come about
and on how the shape of the polytope near any point can be read off
from infinitesimal data on the manifold.  It also applies to some
interesting classes of noncompact or singular Hamiltonian spaces, such
as cotangent bundles and complex affine varieties.
\end{abstract}


\maketitle

\tableofcontents

\section{Introduction}

Let $K$ be a compact Lie group acting smoothly on a compact symplectic
manifold $M$ and suppose there exists a moment(um) map for the action.
This map has a host of interesting properties, one of the most
important of which is the fact that the intersection of its image with
any Weyl chamber is a convex polytope, referred to as the
\emph{momentum polytope\/} of $M$.  This theorem, which is due to
Kirwan, has a long history, part of which I now summarize.  Kostant
proved a convexity theorem for torus actions on conjugacy classes and
flag manifolds in \cite{ko:on}.  Atiyah in \cite{at:co} and Guillemin
and Sternberg in \cite{gu:co1} dealt with the case of general
Hamiltonian torus actions.  In their paper, Guillemin and Sternberg
further proved a convexity theorem for Hamiltonian actions of
arbitrary compact Lie groups on integral K\"ahler manifolds (or
projective manifolds), which was also proved by Mumford in
\cite{ne:st}.  Kirwan subsequently extended this result to Hamiltonian
actions on arbitrary compact symplectic manifolds in \cite{ki:con}.
Many useful refinements in the projective-algebraic case were made
later by Brion in \cite{br:im}.  See \cite{du:co}, \cite{hi:sy} and
\cite{lu:on} for other results and more references.  See
\cite{he:kae}, \cite{kl:st} and \cite{le:no} for some developments
subsequent to the present paper.

A striking difference between Kirwan's general convexity theorem and
the abelian convexity theorem of Atiyah-Guillemin-Sternberg lies in
the fact that the latter offers far more quantitative information on
the shape of the momentum polytope.  For example, in the abelian case
one knows that the vertices of the polytope are images of fixed points
in $M$, and that the shape of the polytope near a vertex can be read
off from the isotropy action on the tangent space at a corresponding
fixed point.  This follows from a combination of the equivariant
Darboux Theorem and Morse theory applied to the components of the
momentum map.

The goal of this paper is to obtain such information in the nonabelian
case as well.  The main result is Theorem \ref{theorem:convex}, which
is a sharpened version of Kirwan's convexity theorem.  Given a point
$m$ in $M$ mapping to a point $\mu$ in the momentum polytope, it
provides a description of the shape of the polytope near $\mu$ in
terms of the action of the stabilizer of $m$ on polynomials on the
tangent space at $m$.  It also states a necessary criterion for $\mu$
to be a vertex, which generalizes the criterion for the abelian case
referred to above.  Other results include convexity theorems for
actions on affine varieties, Theorem \ref{theorem:affine}, and
cotangent bundles, Theorem \ref{theorem:cotangent}.  Theorem
\ref{theorem:projective} describes the relation between the momentum
cone of an affine variety and the momentum polytopes of its projective
closure and the divisor at infinity.

These results are inspired by Brion's treatment of Kirwan's theorem
for projective varieties.  It came as a surprise to me how well
Brion's algebro-geometric techniques can be adapted to a $C^\infty$
setting essentially without sacrificing any of their power.  The main
reason why this is possible is that every point in $M$ possesses an
invariant neighbourhood that is isomorphic as a Hamiltonian
$K$-manifold to (a germ of) a complex quasi-projective variety.

In the language of the orbit method, the momentum polytope of $M$ is
the ``classical'' analogue of the set of highest weights of the
unitary irreducible representations occurring in the ``quantization''
of $M$.  There is also a classical analogue of the space of
highest-weight vectors.  This will be the subject of a forthcoming
paper.

The paper is organized as follows.  Section \ref{section:preliminary}
is a review of some basic facts concerning representations and
momentum maps.  Section \ref{section:semistable} is a review of the
convexity theorem for complex projective varieties, where I have
presented the argument in such a manner that it can be applied to
noncompact varieties.  In Sections \ref{section:affine} and
\ref{section:stein} I prove convexity theorems for complex affine and
Stein varieties.  In Section \ref{section:convex} I apply these
results to prove local convexity properties of arbitrary momentum
maps, whence I derive the main result, Theorem \ref{theorem:convex}.
The local description of the momentum polytope given by this theorem,
although explicit, is unwieldy in practice, and one often has to
revert to ad hoc methods to calculate momentum polytopes.  In Section
\ref{section:examples} I illustrate this in a number of examples, such
as actions on cotangent bundles and projective spaces.

I thank Sheldon Xu-Dong Chang, Yael Karshon and Eugene Lerman for
their help and encouragement.  I am grateful to Laurent Laeng, Domingo
Luna and the referee for correcting a number of errors.

\section{Preliminaries}\label{section:preliminary}

In this section I introduce notation and review basic material to be
used later.

\subsection{Groups, representations}\label{subsection:group}

Throughout this paper $K$ will be a compact connected Lie group with a
fixed maximal torus $T$.  The complexification of $K$ is denoted by
$G$ and the complexification of $T$ by $H$.  Let us fix a Borel
subgroup $B$ of $G$ containing $H$.  Its unipotent radical $(B,B)$ is
denoted by $N$, and the corresponding positive Weyl chamber in $\lie
t^*$ by $\tplus$.  The lattice $\ker(\exp|_{\lie t})$ is denoted by
$\Lambda$.  Its dual lattice
$$
\Lambda^*=\Hom_{\Z}(\Lambda,\Z)\subset\lie t^*
$$ 
is the lattice of (real) weights and $\lplus=\Lambda^*\cap\tplus$ is
the monoid of dominant weights.  To a real weight $\lambda$
corresponds a character $\zeta_\lambda$ of $T$ defined by
$\zeta_\lambda(\exp\xi)=\exp(2\pi\sq\,\langle\lambda,\xi\rangle)$ for
$\xi\in\lie t$.

The complex reductive group $G=K\co$ has a unique complex affine
structure.  Let $R=\C[G]^N$ be the algebra of polynomial functions on
$G$ which are $N$-invariant on the right, that is to say, $f\in R$ if
$f(gn)=f(g)$ for all $g\in G$ and $n\in N$.  Then $G$ acts on $R$ by
left multiplication and, since $H$ normalizes $N$, $H$ acts on $R$ by
right multiplication.  Under the right $H$-action, $R$ has a weight
space decomposition
\begin{equation}\label{equation:ring}
R=\bigoplus_{\lambda\in\lplus}R_\lambda,
\end{equation}
and it follows from the Borel-Weil Theorem that $R_\lambda$ is an
irreducible $G$-module with highest weight $\lambda$.  (See
\cite{kr:ge}, Ch.\ III.)  This implies that the algebra $R$ is of
finite type and, hence, that the scheme $G\qu N=\spec R$ ``is'' an
affine variety.  (I shall not distinguish between an affine variety
and the scheme associated to it.)

Let $\lie W=N(T)/T$ be the Weyl group of $(K,T)$ and let $w_0$ be the
longest Weyl group element.  Define an involution $*\colon\lie
t\to\lie t$ by $\mu^*=-w_0\mu$.  The complex-linear extension of $*$
to $\lie t\co$ and the dual map on $(\lie t^*)\co$ will also be
denoted by $*$.  It is well-known that $*$ leaves the set of dominant
weights invariant and that for all $\lambda\in\lplus$ the
representation $R_{\lambda^*}$ is isomorphic to $R_\lambda^*$, the
contragredient representation of $R_\lambda$.

\subsection{Cones, polytopes}\label{subsection:cone}

Let $E$ be a finite-dimensional vector space over $\Q$ and let $S$ be
a subset of $E$.  The \emph{$\Q$-convex hull\/} of $S$ is the smallest
convex subset of $E$ containing $S$ and is denoted by $\hull_{\Q}S$.
The \emph{convex hull\/} of $S$ is the smallest convex subset of
$E\otimes\R$ containing $S$ and is denoted by $\hull S$.  A subset of
$E$ (resp.\ $E\otimes\R$) is called a \emph{cone\/} if it is invariant
under multiplication by nonnegative rational (resp.\ real) scalars.
The \emph{convex $\Q$-cone\/} spanned by $S$ is the smallest convex
cone in $E$ containing $S$ and is denoted by $\cone_{\Q}S$.  In other
words, $\cone_{\Q}S=\Q_{\geq0}\cdot\hull_{\Q}S$.  The \emph{convex
cone\/} spanned by $S$ is the smallest convex cone in $E\otimes\R$
containing $S$ and is denoted by $\cone S$.  That is, $\cone
S=\R_{\geq0}\cdot\hull S$.  If $S$ is finite, $\hull S$ and $\cone S$
are called a \emph{rational convex polytope}, resp.\ a \emph{rational
convex polyhedral cone\/} in $E\otimes\R$.  A cone is called
\emph{proper\/} if it does not contain any linear subspaces (apart
from $\{0\}$).

\subsection{Hamiltonian actions}\label{subsection:hamilton}

Let $(M,\omega)$ be a symplectic manifold with a $K$-action defined by
a smooth map $\tau\colon K\times M\to M$.  The action $\tau$ is called
\emph{Hamiltonian\/} if there exists a \emph{momentum map}, that is, a
map $\Phi\colon M\to\lie k^*$ with the property that $d\Phi^\xi
=\iota(\xi_M)\omega$ for all $\xi\in\lie k$.  Here $\xi_M$ denotes the
vector field on $M$ induced by $\xi$, and $\Phi^\xi$ is the function
defined by $\Phi^\xi(m)=\bigl(\Phi(m)\bigr)(\xi)$.  We may, and will,
assume $\Phi$ to be $K$-equivariant with respect to the coadjoint
action on $\lie k^*$.  The quadruple $(M,\omega,\tau,\Phi)$ is called
a \emph{Hamiltonian $K$-manifold}.  (See e.g.\ \cite{gu:sy}.)

If $Y$ is any subset of $M$, we denote the restriction of $\Phi$ to
$Y$ by $\Phi_Y$.  The \emph{momentum set\/} $\Delta(Y)$ of $Y$ is
defined by
$$
\Delta(Y)=\Phi(Y)\cap\tplus.
$$
In ``good'' cases, $\Delta(Y)$ is known to be a convex cone or
polytope (\cite{gu:co1}, \cite{ki:con}), and is then called the
\emph{momentum cone}, resp.\ \emph{polytope\/} of $Y$.

On every Hamiltonian $K$-manifold $(M,\omega,\tau,\Phi)$ there exists
an almost-complex structure $J$ that is \emph{compatible\/} with the
Hamiltonian action, that is to say, $J\colon TM\to TM$ is a symplectic
map, the symmetric bilinear form $\omega(\cdot,J\cdot)$ is positive
definite, and $J$ is $K$-equivariant.  If in addition $J$ is
integrable, then $M$ is a K\"ahler manifold with $K$-invariant metric
$ds^2=\omega(\cdot,J\cdot)-\sq\,\omega(\cdot,\cdot)$, and $K$ acts
holomorphically.

\begin{example}[\cite{gu:co1}]\label{example:linear}
Let $(V,\omega_V)$ be a symplectic vector space and assume $K$ acts on
$V$ by linear symplectic transformations.  This action is Hamiltonian;
a momentum map is given by the quadratic map
\begin{equation}\label{equation:quadratic1}
\Phi_V^\xi(v)=\frac1{2}\omega_V\bigl(\xi v,v\bigr),
\end{equation}
where $\xi v$ denotes the image of $v\in V$ under $\xi\in\lie k$
(viewed as a linear operator on $V$).  Choose a $K$-invariant
$\omega_V$-compatible complex structure $J$ on $V$.  Let
$\langle\cdot,\cdot\rangle$ be the Hermitian inner product whose
imaginary part is equal to $-\omega_V$.  Then the momentum map can also
be written as
\begin{equation}\label{equation:quadratic2}
\Phi_V^\xi(v)=\frac\sq{2}\langle\xi v,v\rangle.
\end{equation}
Now suppose $K=T$ is a torus.  Then $V$ is an orthogonal direct sum of
weight spaces, $V=\bigoplus_{\nu\in\Lambda^*}V_\nu$.  If $V_\nu\neq0$,
then $\nu$ is called a weight of the symplectic action of $T$ on $V$.
The weight space decomposition depends on the choice of the complex
structure, but the weights do not.  (This follows from the fact that
any two $K$-invariant compatible complex structures on $V$ are
conjugate by a $K$-equivariant linear symplectic map.)  If $v$ is a
vector of weight $\nu$, then $\xi v=2\pi\sq\,\nu(\xi)v$, so
$\Phi(v)=-\pi\|v\|^2\nu$ by \eqref{equation:quadratic2}, and therefore
\begin{equation}\label{equation:weights}
\Delta(V)=-\cone\{\nu_1,\dots,\nu_l\},
\end{equation}
where $\nu_1,\dots$, $\nu_l$ are the (real) weights of $V$.
\end{example}

\begin{example}[\cite{ki:coh},\cite{ne:st},\cite{ar:co}] 
\label{example:projective}
Let $V$, $J$ and $\langle\cdot,\cdot\rangle$ be as in the previous
example, and let $\P V$ be the space of complex lines in $V$.  The
natural $K$-action on $\P V$ leaves the Fubini-Study symplectic form
invariant, and with the volume of $\P V$ normalized to 1, a
momentum map is given by
\begin{equation}\label{equation:projective1}
\Phi_{\P V}([v])=\frac{\Phi_V(v)}{\pi\|v\|^2}
=\const\frac{\langle\xi v,v\rangle}{\|v\|^2},
\end{equation}
where $[v]$ denotes the line through $v$.  Consequently, if $K$ is a
torus and $v$ is a weight vector in $V$ with weight $\nu$, then
$\Phi_{\P V}([v])=-\nu$.  Hence,
\begin{equation}\label{equation:weyl}
\Delta(\P V)=-\hull\{\nu_1,\dots,\nu_l\},
\end{equation}
where $\nu_1,\dots$, $\nu_l$ are the weights of $V$.  See further
Section \ref{subsection:projective}.
\end{example}

\subsection{Coadjoint orbits}\label{subsection:orbit}

For every $\mu$ in $\tplus$ the coadjoint orbit $K\mu$ with its
Kirillov-Kostant symplectic form $\omega_\mu$ is a Hamiltonian
$K$-manifold.  The momentum map is simply the inclusion
$\iota_\mu\colon K\mu\to\lie k^*$.  (See \cite{ko:qu}, \cite{gu:sy}.)
Let $P_\mu$ be the parabolic subgroup $(K_\mu)\co N$ of $G$.  The
$K$-equivariant diffeomorphism $G/P_\mu\to K\mu$ sending the coset
$1P_\mu$ to the vector $\mu$ provides $K\mu$ with a complex structure
with respect to which $\omega_\mu$ is K\"ahler.

Now let $\mu$ be a dominant weight.  Then the cohomology class of the
form $\omega_\mu$ on $K\mu$ is integral and because $K\mu$ is compact,
there exists a Hermitian holomorphic line bundle $\O_\mu$ on $K\mu$
with curvature $-2\pi\sq\,\omega_\mu$.  The pullback of $\O_\mu$ to
$G/P_\mu$ is just the homogeneous line bundle $G\times^{P_\mu}\C$,
where $P_\mu$ acts on $\C$ by the character $\mu$.

Let $\mu_1$ and $\mu_2$ be two points in $\tplus$ and let
$\mu=\mu_1+\mu_2$.  Then $K_\mu=K_{\mu_1}\cap K_{\mu_2}$ and
$P_\mu=P_{\mu_1}\cap P_{\mu_2}$, so we have canonical holomorphically
locally trivial fibrations $\pi_i\colon K\mu\to K\mu_i$.  It is not
hard to see that $\omega_\mu =\pi_1^*\omega_{\mu_1}
+\pi_2^*\omega_{\mu_2}$.  If $\mu$, $\mu_1$ and $\mu_2$ are dominant,
the holomorphic line bundle $\O_\mu$ is isomorphic to
$\pi_1^*\O_{\mu_1}\otimes\pi_2^*\O_{\mu_2}$.  Let me summarize this in
a commutative diagram:
$$
\xymatrix{
&
{\pi_1^*\O_{\mu_1}\oplus\pi_2^*\O_{\mu_2}}\ar[dr]\ar[d]^\otimes\ar[dl]
\\
{\O_{\mu_1}}\ar[d] & {\O_\mu}\ar[d] & {\O_{\mu_2}}\ar[d] \\
K\mu_1 & K\mu\ar[l]_{\pi_1}\ar[r]^{\pi_2} & K\mu_2.
}
$$

\subsection{Gradient flows, semistability}\label{subsection:flow}

Let $X$ be a smooth connected Riemannian manifold and let $f$ be a
function on $X$ with the property that at every point of $X$ there
exists a system of local coordinates in which $f$ is real-analytic.
Let $\lie F(t,\cdot)$ be the gradient flow of $-f$.  Assume that the
path of steepest descent $\lie F(t,x)$ through every point $x$ is
contained in a compact set.  Then the flow is defined for all
$t\geq0$.  Moreover, by results of {\L}ojasiewicz \cite{lo:pr} and
Simon \cite{si:as}, the limit $x_\infty =\lim_{t\to\infty}\lie F(t,x)$
exists for all $x$.  Let $a$ be a critical level of $f$ and let
$S_a=\{x\in X:f(x_\infty)=a\}$ be the stable set of $a$.  Then $S_a$
is a locally closed subset of $X$ and that $x\mapsto x_\infty$ is a
continuous retraction from $S_a$ onto $f\inv(a)$.  The decomposition
$X=\coprod_aS_a$ is called the \emph{Morse decomposition\/} of $X$
with respect to $f$ (even if $f$ is not a Morse function).

(A note on the literature: {\L}ojasiewicz' paper \cite{lo:pr} does not
contain a complete proof of these assertions.  He explains part of the
requisite estimates in \cite{lo:en}.  Simon gives a fuller account of
the retraction argument in \cite{si:as}, while also generalizing it to
an infinite-dimensional situation.  Apparently independently of both
\cite{lo:pr} and \cite{si:as}, Neeman \cite{ne:to} and Schwarz
\cite{sc:to} rederive the retraction theorem for certain flows on
vector spaces from the inequalities in \cite{lo:en}.  Their arguments
can easily be generalized to prove the above statements.)

As an example, let $(X,\omega,\tau,\Phi)$ be a connected Hamiltonian
$K$-manifold equipped with a compatible almost-complex structure $J$.
Let $(\cdot,\cdot)$ be a $K$-invariant inner product on $\lie k$ and
let $\lvert\cdot\rvert$ be the associated norm.  I use the same
symbols to denote the corresponding inner product and norm on the dual
$\lie k^*$.  Put $f=|\Phi|^2$.  It follows from the local model for
Hamiltonian actions (see Section \ref{section:convex}) that this
function is real-analytic in suitable local coordinates.  Note that
since $f$ is $K$-invariant, the flow $\lie F(t,\cdot)$ is
$K$-equivariant.  Let us assume the momentum map to be
\emph{admissible\/} in the sense that for every $x\in X$ the path of
steepest descent $\lie F(t,x)$ is contained in a compact set.  The set
$S_0$ is called the set of (\emph{analytically\upn) semistable\/}
points and is denoted by $X\sst(\Phi)$.  So $X\sst(\Phi)$ is nonempty
if and only if $0\in\Phi(X)$.  Kirwan has shown in \cite{ki:coh} that
$S_a$ is a submanifold of even codimension for every critical level
$a$ (and if $J$ is integrable, then $S_a$ is a complex submanifold).
Now assume that the Morse decomposition $X=\coprod_aS_a$ is locally
finite.  (This is for instance the case if for every $a$ there are
only finitely many critical levels below $a$).  Then, if nonempty,
$X\sst(\Phi)$ is open, connected and dense.

\section{Semistability and convexity}\label{section:semistable}

In this section I review the convexity theorem for K\"ahler manifolds
due to Guillemin and Sternberg \cite{gu:co1} and Mumford
\cite{ne:st}.  Guillemin and Sternberg have pointed out in
\cite{gu:co2} that semistability and convexity are closely
related.  This idea goes back to Heckman's paper \cite{he:pr}, and can
be formulated as follows.

\begin{proposition}\label{proposition:semistable}
Let $(Y_i,\sigma_i,\Psi_i)$ be compact Hamiltonian $K$-manifolds with
compatible \upn(integrable\upn) complex structures $J_i$\upn, where
$i=1$\upn, $2,\dots$\upn, $k$.  Assume that the cohomology classes of
the $\sigma_i$ are integral.  Let $Y$ be a compact complex
$K$-manifold and let $p_i\colon Y\to Y_i$ be $K$-equivariant
surjective holomorphic maps.  Let $a_1$\upn, $a_2,\dots$\upn, $a_k$ be
nonnegative numbers\upn, and let $\sigma=\sum_ia_ip_i^*\sigma_i$ and
$\Psi=\sum_ia_ip_i^*\Psi_i$.  Assume $\sigma$ is a K\"ahler form on
$Y$.  Then $\bigcap_i p_i\inv\bigl(Y\sst(\Psi_i)\bigr)$ is contained
in $Y\sst(\Psi)$.  Hence\upn, if $0\in\Psi_i(Y_i)$ for every $i$\upn,
then $0\in\Psi(Y)$.
\end{proposition}

By the equivariant version of Kodaira's Embedding Theorem, the
manifolds $Y_i$ and $Y$ are of course biholomorphically equivalent
(but not necessarily isometric) to projective manifolds with linear
$G$-actions.  The proof is a straightforward adaptation of the
techniques of \cite{gu:co1} and \cite{ne:st}.  With a view to later
applications I supply an argument which can easily be made to work for
noncompact manifolds.

\begin{proof}
Note that $\Psi$ is a momentum map for the $K$-action on $Y$ with
respect to the symplectic form $\sigma$.  Also, since $Y$ is compact,
$\Psi$ is admissible in the sense of Section \ref{subsection:flow}.
For clarity I will first handle the case where $Y_i=Y$ and
$p_i=\id_Y$.  So we are given K\"ahler forms $\sigma_i$ on $Y$ and we
wish to show $\bigcap_i Y\sst(\Psi_i)\subset Y\sst(\Psi)$.

For $i=1$, $2,\dots$, $k$, let $L_i$ be a Hermitian holomorphic line
bundle on $Y$ with curvature form $-2\pi\sq\,\sigma_i$.  (These exist
because $Y$ is compact.)  Let $n_i$ be a positive integer and let
$s_i\in\Gamma(Y,L_i^{n_i})^K$, where $\Gamma$ stands for holomorphic
sections.  Let $\langle s_i,s_i\rangle$ denote the length squared of
$s_i$ with respect to the fibre metric on $L_i^{n_i}$.  It follows
from the invariance of the $s_i$ that for every $\xi\in\lie k$ we have
$\Lie_{J\xi_Y}\langle s_i,s_i\rangle=-4\pi n_i\Psi_i^\xi\langle
s_i,s_i\rangle$, and hence
\begin{equation}\label{equation:length}
\Lie_{J\xi_Y}\langle s_i,s_i\rangle^{a_i/n_i}=-4\pi
a_i\Psi_i^\xi\langle s_i,s_i\rangle^{a_i/n_i}.
\end{equation}
Here $\Lie$ stands for the Lie derivative.  Let $L$ be the line bundle
$\bigotimes_iL_i$ with the product Hermitian metric and let $\lie
F(t,\cdot)$ be the flow of the vector field $-\grad|\Psi|^2$.  From
the elementary fact that
\begin{equation}\label{equation:gradient}
J\xi_Y=\grad\Psi^\xi
\end{equation}
one easily deduces that
\begin{equation}\label{equation:yangmills}
\grad|\Psi(y)|^2=2J\Psi(y)^\flat_{Y,y},
\end{equation}
where $\flat\colon\lie k^*\to\lie k$ is the linear isomorphism defined
by the inner product, and $\Psi(y)^\flat_Y$ is the vector field on $Y$
induced by $\Psi(y)^\flat$.  (See \cite{ki:coh}.)  Put
$s=s_1^{n/n_1}\otimes s_2^{n/n_2}\otimes\cdots\otimes s_k^{n/n_k}$,
where $n=n_1n_2\cdots n_k$.  Then $s\in\Gamma(Y,L^n)^K$.  Consider the
function $u=\langle s_1,s_1\rangle^{a_1/n_1}\langle
s_2,s_2\rangle^{a_2/n_2}\cdots\langle s_k,s_k\rangle^{a_k/n_k}$.  By
\eqref{equation:length},
%
$$
\Lie_{J\xi_Y}u
=-4\pi(a_1\Psi_1^\xi+a_2\Psi_2^\xi+\cdots+a_k\Psi_k^\xi)u
=-4\pi\Psi^\xi u.
$$
Using this and \eqref{equation:yangmills} we find that the derivative
of $u$ along a trajectory $\lie F(t,y)$ is equal to
\begin{multline}\label{equation:maximum}
\frac{d}{dt}u\bigl(\lie F(t,y)\bigr)=-du\bigl(\grad|\Psi(\lie
F(t,y))|^2\bigr) =-2du\Bigl(J\Psi\bigl(\lie F(t,y)\bigr)^\flat_{Y,\lie
F(t,y)}\Bigr) \\=8\pi\bigl\langle\Psi\bigl(\lie F(t,y)\bigr),
\Psi\bigl(\lie F(t,y)\bigr)\bigr\rangle u\bigl(\lie F(t,y)\bigr)
=8\pi\bigl|\Psi\bigl(\lie F(t,y)\bigr)\bigr|^2u\bigl(\lie
F(t,y)\bigr)\geq0.
\end{multline}
Now suppose $y\in Y$ is semistable with respect to all $\sigma_i$.
Kirwan \cite{ki:coh} and Ness \cite{ne:st} observed that for a
projective manifold with the Fubini-Study metric analytic
semistability is equivalent to semistability in Mumford's sense.  This
is true in general for a compact complex manifold with an integral
K\"ahler metric; see \cite{sj:ho}.  This means we can find positive
integers $n_i$ and invariant global holomorphic sections $s_i$ of
$L_i^{n_i}$ such that $s_i(y)\neq0$.  Then
$s(y)=s_1^{n/n_1}(y)\otimes\cdots\otimes s_k^{n/n_k}(y)\neq0$, so
$u(y)>0$.  Put $y_\infty=\lim_{t\to\infty}\lie F(t,y)$.  By
\eqref{equation:maximum}, $u\bigl(\lie F(t,y)\bigr)$ is increasing
along the path $\lie F(t,y)$, so $u(y_\infty)>0$.  On the other hand,
$du\bigl(\lie F(t,y)\bigr)\big/dt$ tends to zero as $t$ tends to
infinity, so from \eqref{equation:maximum} we get
$|\Psi(y_\infty)|^2u(y_\infty)=0$, and therefore
$|\Psi(y_\infty)|^2=0$.  In other words, $y$ is semistable for
$\sigma$.  We have shown $\bigcap_iY\sst(\Psi_i)\subset Y\sst(\Psi)$.

Suppose now that $0\in\Psi_i(Y_i)$ for all $i$.  Then the sets
$Y\sst(\Psi_i)$ are nonempty for all $i$, and are therefore open and
dense.  It follows that their intersection is nonempty.  If
$y\in\bigcap_iY\sst(\Psi_i)$, then, by the first part of the proof,
$\Psi(y_\infty)=0$, that is, $0\in\Psi(Y)$.

In the general case the argument is almost exactly the same.  The
difference is that one considers the line bundle
$L=\bigotimes_ip_i^*L_i$ on $Y$ and the section
$s=p_1^*s_1^{n/n_1}\otimes\cdots\otimes p_k^*s_k^{n/n_k}$ of $L^n$.
Further, the assumptions on the $p_i$ guarantee that the pre-image
$p_i\inv(S)$ of a complex-analytic subset $S\subset Y_i$ of positive
codimension is a complex-analytic subset of positive codimension of
$Y$.  This implies that if $Y_i\sst\neq\emptyset$ for all $i$, then
$\bigcap_ip_i\inv\bigl(Y\sst(\Psi_i)\bigr)\neq\emptyset$.
\end{proof}

For noncompact $Y$ the flow of $-\grad|\Psi|^2$ may not be defined for
all time or its trajectories may fail to converge, and the equivalence
between analytic and algebraic semistability can break down.  But the
following qualified statement is still true.  The proof is almost word
for word the same.

\begin{proposition}\label{proposition:noncompact}
Let $(Y_i,\sigma_i,\Psi_i)$ be Hamiltonian $K$-manifolds endowed with
compatible complex structures $J_i$\upn, where $i=1$\upn,
$2,\dots$\upn, $k$.  Assume there exist $K$-equivariant Hermitian
holomorphic line bundles $L_i$ on $Y_i$ with curvature forms
$-2\pi\sq\,\sigma_i$ for all $i$.  Let $Y$ be a complex $K$-manifold
and let $p_i\colon Y\to Y_i$ be $K$-equivariant surjective holomorphic
maps.  Let $a_1$\upn, $a_2,\dots$\upn, $a_k$ be nonnegative
numbers\upn, and let $\sigma=\sum_ia_ip_i^*\sigma_i$ and
$\Psi=\sum_ia_ip_i^*\Psi_i$.   Assume $\sigma$ is a K\"ahler form on
$Y$ and that the momentum map $\Psi$ is admissible in the sense of
Section \ref{subsection:flow}.  If for every $i$ there exist a positive
integer $n_i$ and a nonzero $K$-invariant holomorphic section of
$L_i^{n_i}$\upn, then $0\in\Psi(Y)$.
\qed
\end{proposition}
 
\begin{remark}\label{remark:singular}
Suppose that $Z\subset Y$ and $Z_i\subset Y_i$ are irreducible
$K$-stable locally closed analytic subvarieties, and that for all $i$
the restriction of the map $p_i$ to $Z$ is a surjective map $Z\to
Z_i$.  Assume that for every $z$ in $Z$ the path $\lie F(t,z)$ and its
limit $z_\infty$ are contained in $Z$.  Also assume that the sections
$s_i$ restrict to nonzero sections on $Z_i$.  Then $0\in\Psi_i(Z_i)$
for all $i$ implies $0\in\Psi(Z)$.  Exactly the same proof works.
\end{remark}

Here is an application of Proposition \ref{proposition:noncompact},
where the notation is as in Section \ref{subsection:orbit}.  Let
$(X,\omega,\Phi)$ be a Hamiltonian $K$-manifold, not necessarily
compact, with a compatible complex structure $J$.  Suppose that there
exists a $K$-equivariant Hermitian holomorphic line bundle $L$ on $X$
with curvature form $-2\pi\sq\,\omega$.  For $i=1$, $2,\dots$, $k$,
let $\mu_i$ be a dominant weight and let $Y_i$ be the manifold
$X\times K\mu_i^*$ with symplectic form
$\sigma_i=\omega+\omega_{\mu_i^*}$.  (Recall that $*$ is defined by
$\nu^*=-w_0\nu$.)  Then the $K$-action on $(Y_i,\sigma_i)$ is
Hamiltonian with momentum map $\Psi_i=\Phi+\iota_{\mu_i^*}$.  Let
$L_i$ be the Hermitian line bundle $L\otimes\O_{\mu_i^*}$ on $Y_i$.
Let $a_i$ be arbitrary positive numbers, let $\mu=\sum_ia_i\mu_i$, and
let $Y$ be the $K$-manifold $X\times K\mu^*$.  Consider the fibrations
$p_i\colon Y\to Y_i$ induced by the fibrations of coadjoint orbits
$K\mu^*\to K(a_i\mu_i)\xrightarrow{q}K\mu_i^*$, where $q$ is the
equivariant diffeomorphism sending $a_i\mu_i$ to $\mu_i$.  Since the
$a_i$ are positive, the form $\sigma=\sum_ia_ip_i^*\sigma_i$ is a
K\"ahler form on $Y$, and the action on $Y$ is Hamiltonian with
momentum map $\Psi=\sum_ia_ip_i^*\Psi_i$.  Let us assume that
\begin{gather} 
\label{assumption:nonvanishing} \text{for all $i$ there exist $n_i>0$
and nonvanishing sections $s_i\in\Gamma(Y_i,L_i^{n_i})^K$;}\\
\label{assumption:admissible} \text{for all $a_i\geq0$ the momentum
map $\Psi$ is admissible.}
\end{gather}
(Assumption \ref{assumption:admissible} holds e.g.\ when $\Phi$ is
proper.)  Then by Proposition \ref{proposition:noncompact} there exists
a point $(x,k\mu^*)$ in $Y=X\times K\mu^*$ with $\Psi(x,k\mu^*)=0$,
that is, $(a_1+\cdots+a_k)\Phi(x)=
kw_0(a_1\mu_1+\cdots+a_k\mu_k)$.  This shows that
$$
\frac{a_1\mu_1+a_2\mu_2+\cdots+a_k\mu_k}{a_1+a_2+\cdots+a_k}
\in\Phi(X)\cap\tplus=\Delta(X)
$$
for all $a_i>0$.  The same is obviously true if some of the $a_i$ are
$0$ (just replace the $\lambda_i$ by the subset consisting of those
$\mu_i$ for which $a_i\neq0$), so we have proved:

\begin{proposition}\label{proposition:tortuous}
Under the assumptions \eqref{assumption:nonvanishing} and
\eqref{assumption:admissible} the convex hull of $\mu_1$\upn,
$\mu_2,\dots$\upn, $\mu_k$ is contained in $\Delta(X)$.
\qed
\end{proposition}

\begin{remark}\label{remark:saturated}
By Remark \ref{remark:singular}, a similar result holds for an
irreducible $K$-stable locally closed analytic subvariety $Z$ of $X$.
Namely, assume that the gradient flow of the function $-|\Psi|^2$
preserves the subvariety $Z\times K\mu^*\subset X\times K\mu^*$ and
that the forward trajectories converge to points in $Z\times K\mu^*$.
(This assumption is satisfied e.g.\ if the $K$-action on $X$ extends
to a holomorphic $G$-action and if $Z$ has the property that
$\overline{Gz}\subset Z$ whenever $z\in Z$.)  Assume further that for
all $i$ the sections $s_i$ restrict to nonzero sections on $Z\times
K\mu_i$.  Then $\hull\{\mu_1,\mu_2,\dots,\mu_k\}\subset\Delta(Z)$.
\end{remark}

Now consider the graded $G$-algebra
$\eu A=\bigoplus_{n\geq0}\Gamma(X,L^n)$.

\begin{definition}[Brion]\label{definition:brion}
The \emph{highest-weight set\/} of $X$ is the subset $\hw(X)$ of the
$\Q$-vector space $\Lambda^*\otimes\Q$ consisting of all $\lambda/n$
with the property that $\lambda$ is a dominant weight and $n$ a
positive integer such that the irreducible representation $R_\lambda$
occurs in the degree-$n$ piece $\eu A_n$.  In other words,
$\lambda/n\in\hw(X)$ if and only if $n>0$ and $\lambda$ occurs as a
weight of the $T$-action on the degree-$n$ piece of the ring $\eu
A^N$.
\end{definition}

Let me recall a few basic facts concerning highest-weight sets.  See
\cite{br:im} for details.  Note first that $\hw(X)$ is contained in
$\tplus$.  The fact that $\eu A$ has no zero divisors implies that
$\hw(X)$ is a convex subset of the $\Q$-vector space
$\Lambda^*\otimes\Q$.  Now suppose $X$ is compact.  Then $\eu A$ is of
finite type.  A result of Luna and Vust says that for every
$G$-algebra $\eu A$ of finite type (not necessarily graded)
\begin{equation}\label{equation:unipotent}
\eu A^N\cong(R\otimes \eu A)^G,
\end{equation}
where $R$ is as in \eqref{equation:ring}.  This implies that $\eu A^N$
is of finite type.  (See e.g.\ \cite{kr:ge}.)  It follows from this
that $\hw(X)$ is the convex hull of a finite number of points in
$\Lambda^*\otimes\Q$.  Moreover, $\Delta(X)$ is closed.  It is now
easy to deduce the following result from Proposition
\ref{proposition:tortuous}.

\begin{theorem}[\cite{gu:co1},\cite{ne:st},\cite{br:im}]
\label{theorem:kaehler} 
If $X$ is a compact integral Hamiltonian $K$-manifold with a
compatible K\"ahler structure\upn, $\Delta(X)$ is equal to the closure
of $\hw(X)$ in $\tplus$ and is therefore a rational convex polytope.
\qed
\end{theorem}

By Remark \ref{remark:saturated}, Theorem \ref{theorem:kaehler} also
holds if we replace $X$ with a $G$-stable irreducible closed analytic
subvariety.

This result applies in particular to a $G$-stable irreducible closed
subvariety $X$ of $\P V$, the projective space of a $G$-module $V$,
equipped with the Fubini-Study symplectic form.  But now consider a
subvariety $X$ of $\P V$ that is not necessarily irreducible or even
reduced.  What is the connection between $\Delta(X)$ and $\hw(X)$?
Let $X\red$ denote the reduced variety associated to $X$ and let
$X_1$, $X_2,\dots$, $X_l$ be its irreducible components (each endowed
with the reduced induced structure).  Let $Z$ and $Z_i$ be the affine
cones of $X\red$, resp.\ $X_i$, and let $I$, resp.\ $I_i$, be their
homogeneous ideals.  Then $I=\bigcap_iI_i$ and for each $i$ we have an
exact sequence
$$
\xymatrix{I_i/I\ar@{((->}[r] & {\C}[Z]\ar@{->>}[r] & {\C}[Z_i]}.
$$
It follows easily from this that $\hw(X\red)=\bigcup_i\hw(X_i)$.  It
is evident that $\Delta(X\red)=\bigcup_i\Delta(X_i)$, so applying
Theorem \ref{theorem:kaehler} to each of the $X_i$ we obtain the
following result.

\begin{addendum}[\cite{la:hi}]\label{addendum:unreduced}
If $X$ is a \upn(not necessarily irreducible or reduced\/\upn)
subvariety of $\P V$\upn, then $\Delta(X\red)$ is equal to the closure
of $\hw(X\red)$ in $\tplus$ and is therefore a union of finitely many
rational convex polytopes.
\qed
\end{addendum}

The following example shows that $\hw(X)$ is in general not a union of
finitely many $\Q$-convex sets and that it is not necessarily equal to
$\hw(X\red)$.

\begin{example}\label{example:unreduced}
Fix a positive integer $n$.  Let $\eu A$ be the algebra
$\C[x,y]/(x^n)$ and let $X=\proj \eu A$, considered as a closed
subscheme of $\P^1$.  For any weight $\lambda$ of $K=S^1$ define an
action of $K$ on $\C^2$ by $g(x,y)=(\zeta_{-\lambda}(g)x,y)$, where
$\zeta_{-\lambda}$ is the character defined by
$\zeta_{-\lambda}(\exp\xi)=
\exp\bigl(-2\pi\sq\langle\lambda,\xi\rangle\bigr)$.  According to
\eqref{equation:weyl}, $\Delta(\P^1)$ is equal to the interval between
$0$ and $\lambda$.  Note that $X\red$ is the point with homogeneous
coordinates $[0,1]$, so
$$
\hw(X\red)=\Delta(X\red)=\{0\}.
$$
Write $\bar f$ for $f+(x^n)\in \eu A$; then $\bar{x}^k\bar{y}^l$ has
degree $k+l$ and weight $k\lambda$ for $k=0$, $1,\dots$, $n-1$ and
$l\ge0$, so that
$$
\hw(X)=\biggl\{\,\frac{k\lambda}{k+l}:\text{$k$, $l\in\N$, $k<n$,
$kl\ne0$}\,\biggr\},
$$
which cannot be written as a union of finitely many $\Q$-convex sets.
\end{example}

Nevertheless, it is always true that the intersection of $\hw(X)$ with
a rational line in $\lie t^*$ is a bounded set which contains both its
endpoints.

\begin{addendum}\label{addendum:sup}
For every \upn(not necessarily irreducible or reduced\/\upn)
subvariety $X$ of $\P V$ and every $\nu$ in $\hw(X)$\upn, let
$I_\nu=\{\,q\in\Q:q\nu\in\hw(X)\,\}$.  Then $\inf I_\nu$ and $\sup
I_\nu$ are in $I_\nu$.
\end{addendum}

\begin{proof}
This is similar to the proof of Theorem \ref{theorem:kaehler}.  The
assertion is trivial for $\nu=0$, so let me assume that $\nu\ne0$.
Let $\lie t_1$ be the subalgebra of $\lie t$ annihilated by $\nu$;
then $T_1=\exp\lie t_1$ is a subtorus of $T$ of codimension one.
Denote the quotient circle $T/T_1$ by $T_2$ and identify $\lie t_2^*$,
the kernel of the canonical projection $\lie t^*\to\lie t_1^*$, with
the line $\R\nu\subset\lie t^*$.  Put $\hw_\nu=\hw(X)\cap\R\nu$.

Let $Z$ be the affine cone on $X$ and let $\eu A$ be the graded
algebra $\C[Z]^N$, which as noted above is of finite type.  Since the
maximal torus $T$ normalizes $N$, it acts in a natural way on $\eu A$.
The algebra $\eu B=\eu A^{T_1}$ is likewise of finite type and it
carries a representation of $T_2$.  Let $\lambda_0$ be the (unique)
primitive element of $\Lambda^*$ such that $\lambda_0=p\nu$ for some
positive rational $p$.  Note that $f\in \eu A$ is a weight vector of
weight proportional to $\lambda_0$ if and only if $f$ is in $\eu B$
and is a weight vector for $T_2$.  It follows that $\hw_\nu$ is equal
to the set of $m\lambda_0/n$ such that there exists $f\in\eu B$ of
weight $m\lambda_0$ and degree $n$.  Choose (nonzero) generators
$f_1$, $f_2,\dots$, $f_l$ of $\eu B$ with weights $m_1\lambda_0$,
$m_2\lambda_0,\dots$, $m_l\lambda_0$ and degrees $n_1$, $n_2,\dots$,
$n_l$.  Then $m_j\lambda_0/n_j\in\hw_\nu$, so $m_j/n_j\in I_\nu$ for
$1\le j\le l$.  Moreover, $q\nu\in\hw_\nu$ if and only if there exist
$a_1$, $a_2,\dots$, $a_l$ such that the monomial $\prod_jf_j^{a_j}$ is
nonzero and
$$
q\nu=\frac{\sum_ja_jm_j}{\sum_ja_jn_j}\lambda_0.
$$
In other words, every element of $\hw_\nu$ is a convex combination of
the $m_j\lambda_0/n_j$ (but not every such combination need be in
$\hw_\nu$, because $\eu B$ may have zero divisors).  Now let $r$ be
the minimum of the $m_j/n_j$ and $s$ their maximum.  Then $r$ and $s$
are in $I_\nu$ and $r=\inf I_\nu$ and $s=\sup I_\nu$.
\end{proof}

\section{Affine varieties}\label{section:affine}

In this section $X$ denotes an affine algebraic variety (not
necessarily reduced or irreducible) on which $G$ acts algebraically.
The main results are theorems describing the momentum map image of $X$
with respect to suitable $K$-invariant symplectic forms, Theorems
\ref{theorem:affine} and \ref{theorem:level}.  My main interest is in
smooth varieties, but the proofs turn out to be no harder for general
varieties.  There are some examples at the end of Section
\ref{subsection:highest}.

\subsection{Highest weights and the momentum cone}
\label{subsection:highest}

The natural analogue of Definition \ref{definition:brion} is the
following.

\begin{definition}[Brion]
The \emph{set of highest weights\/} of $X$ is the subset $\hw(X)$ of
$\lplus$ consisting of all dominant weights $\lambda$ such that the
irreducible $G$-representation $R_\lambda$ occurs in the coordinate
ring $\C[X]$.  In other words, $\lambda\in\hw(X)$ if and only if
$\lambda$ is a weight of the $T$-action on the ring $\C[X]^N$.  If $G$
is a torus, we refer to $\hw(X)$ as the \emph{weight set\/} of $X$.
\end{definition}

Note that if $Y$ is another affine $G$-variety, then $\hw(X)=\mathcal
C(Y)$ if and only if the coordinate rings of $X$ and $Y$ contain the
same irreducible $G$-representations (up to multiplicities).  It is
easy to see that if $X$ is irreducible, then $\hw(X)$ is a submonoid
of $\lplus$ (that is, it contains 0 and is invariant under addition).
It follows from \eqref{equation:unipotent} that $\hw(X)$ is finitely
generated as a monoid.

Let me first discuss a few simple tricks for computing highest-weight
sets.  Let $X\qu G=\spec\C[X]^G$ denote the categorical quotient of
$X$, that is, the variety of closed $G$-orbits in $X$, and let
$\pi\colon X\to X\qu G$ be the quotient mapping.  (See \cite{kr:ge} or
\cite{lu:sl}.)  A subset $U$ of $X$ is called \emph{saturated\/} (with
respect to $\pi$) if $\pi\inv\pi(U)=U$, that is, $\overline{Gx}\subset
U$ whenever $x\in U$.  Saturated subsets are evidently $G$-stable.
The following result says that $\hw(X)$ is determined locally in the
Zariski topology, or, more precisely, that it does not change when we
remove from $X$ a divisor defined by an invariant polynomial.

\begin{lemma}\label{lemma:open}
Assume $X$ is irreducible.  Let $Y$ be any saturated affine
Zariski-open subvariety of $X$.  Then $\hw(X)=\hw(Y)$.
\end{lemma}

\begin{proof}
The coordinate ring of $X$ embeds equivariantly into the coordinate
ring of $Y$.  This implies $\hw(X)\subset\hw(Y)$.  For the reverse
inclusion, note that the assumption that $Y$ is saturated implies that
the quotient $Y\qu G$ is an affine open subvariety of $X\qu G$.  Let
$D$ be the complement of $Y\qu G$ in $X\qu G$, let $D_1$, $D_2,\dots$,
$D_k$ be the irreducible components of $D$, and let $f_i\in\C[X\qu G]$
be the polynomial defining $D_i$ for $i=1$, $2,\dots$, $k$.  Put
$f=f_1f_2\cdots f_k$; then the divisor $X-Y=\pi\inv(D)$ is the zero
set of $f$ (viewed as an element of $\C[X]$), and the coordinate ring
of $Y$ is just the localization of $\C[X]$ at $f$.  Define for every
$p$ the linear map $\psi_p\colon\C[X]_p\to\C[Y]$ by $\psi_p(a)=a/f^p$.
Since $X$ is irreducible, $\psi_p$ is injective.  Since $f$ is
invariant, $\psi_p$ is equivariant.  The direct sum of the maps
$\psi_p$ is an equivariant map from $\C[X]$ to $\C[Y]$, which is
surjective, because $\C[Y]$ is the localization of $\C[X]$ at $f$.  It
follows from this that if an irreducible $G$-representation occurs in
$\C[Y]$, then it occurs in $\C[X]$.  This proves
$\hw(Y)\subset\hw(X)$.
\end{proof}

The problem of computing highest-weight sets can in principle be
reduced to torus actions.  The variety $\spec\C[X]^N$ is a categorical
quotient of $X$ by $N$ in the category of affine varieties, and will
be denoted by $X\qu N$.  (It is not always the same as the
set-theoretical quotient of $X$ by $N$ and the natural map $X\to X\qu
N$ is not always surjective.  For instance, the homogeneous space
$G/N$ is a quasi-affine variety and $G\qu N$ is its affine closure.)
Let $\hw(X\qu N)$ be the weight set of the $H$-action on $X\qu N$.  By
the theorem of the highest weight, a weight occurs in $\C[X]^N$ if and
only if it occurs as the highest weight of an irreducible component of
$\C[X]$.  This proves the following lemma.

\begin{lemma}\label{lemma:unipotent}
$\hw(X)=\hw(X\qu N)$.
\qed
\end{lemma}

Furthermore, highest-weight sets are invariant up to denominators
under finite morphisms.

\begin{lemma}\label{lemma:finite}
Let $X$ and $Y$ be affine $G$-varieties and let $\phi\colon X\to Y$ be
a finite surjective $G$-morphism.  Then $\hw(Y)$ is contained in
$\hw(X)$ and $n!\,\hw(X)$ is contained in $\hw(Y)$\upn, where
$n$ is the cardinality of the generic fibre of $\phi$.
\end{lemma}

\begin{proof}
Let $\lie A$ and $\lie B$ be the coordinate rings of $X$, resp.  $Y$.
Then $\lie B$ can be regarded as a subring of $\lie A$ via the
pull-back map $\phi^*$.  This implies $\hw(Y)\subset\hw(X)$.  Now for
the second inclusion.  By Lemma \ref{lemma:unipotent}, it suffices to
show that $n!\,\hw(X\qu N)\subset\hw(Y\qu N)$.  Recall that the
finiteness of $\phi$ means that $\lie A$ is a $\lie B$-module of rank
$n$.  By Satz 1 on p.\ 192 of \cite{kr:ge}, $\lie A^N$ is a $\lie
B^N$-module of rank $\leq n$.  This implies that every element $a$ of
$\lie A^N$ satisfies an equation $P(a)=0$, where $P(t)$ is a monic
polynomial of degree $n$ in $\lie B^N[t]$.  Let $a\in\lie A^N$ be an
element of weight $\lambda\in\Lambda^*$.  I will show that $\lie B^N$
contains an element of weight $k\lambda$ for some $k\leq n$.  There
exist $b_0$, $b_1,\dots$,~$b_{n-1}$ in $\lie B^N$ such that
$a^n+b_{n-1}a^{n-1}+\cdots+b_1a+b_0=0$.  Let $k$ be the largest number
$l$ such that $b_{n-l}\neq0$.  Then
$$
a^k+b_{n-1}a^{k-1}+\cdots+b_{n-k+1}a+b_{n-k}=0.
$$
Because the action of $H$ on $\lie B^N$ is completely reducible, we
may assume all terms in this equation have the same weight.  Then the
weight of $b_{n-k}$ is equal to the weight of $a^k$, which is
$k\lambda$.
\end{proof}

\begin{example}\label{example:galois}
Let $\Gamma$ be a finite group acting on $X$ and suppose that the
actions of $G$ and $\Gamma$ commute.  Let $Y$ be the affine $G$-variety
$X/\Gamma$.  Then the lemma shows that $\hw(Y)\subset\hw(X)$ and
$n!\,\hw(X)\subset\hw(Y)$, where $n$ is the cardinality of
$\Gamma$.
\end{example}

\begin{remark}\label{remark:surjective}
If $\phi\colon X\to Y$ is any surjective $G$-morphism of affine
$G$-varieties, then $\hw(Y)$ is a subset of $\hw(X)$.
\end{remark}

\begin{remark}\label{remark:product}
Suppose that $G$ is the direct product of two reductive subgroups
$G_1$ and $G_2$.  Then the monoid of dominant weights of $G$ is simply
the product of the monoids of dominant weights of $G_1$ and $G_2$, and
the positive Weyl chamber of $G$ is the product of the positive Weyl
chambers of $G_1$ and $G_2$.  Moreover, every irreducible
representation of $G$ is a tensor product of an irreducible
representation of $G_1$ and an irreducible representation of $G_2$.
It follows that for every affine $G_1$-variety $X_1$ and for every
affine $G_2$-variety $X_2$, $\hw(X_1\times X_2)$ is the product
of $\hw(X_1)$ and $\hw(X_2)$.
\end{remark}


Now take any equivariant algebraic closed embedding of $X$ into some
representation space $V$.  Such embeddings always exist; see e.g.\
\cite{kr:ge}.  Let $\langle\cdot,\cdot\rangle$ be a $K$-invariant
Hermitian inner product on $V$ and let $\omega_V$ be the symplectic
form $-\Im\langle\cdot,\cdot\rangle$.  Denote by $\lVert{\cdot}\rVert$
the corresponding norm, and by $\Phi_V$ the momentum map given by
\eqref{equation:quadratic1}.  Now attach a copy of the one-dimensional
trivial representation $\C$ to $V$.  Then the projective space
$\P(V\oplus\C)$ carries a natural $G$-action, and the projective space
$\P V$ can be identified equivariantly with the hyperplane at infinity
in $\P(V\oplus\C)$.  By \eqref{equation:projective1}, the momentum map
on $\P(V\oplus\C)$ is given by
\begin{equation}\label{equation:projective2}
\Phi_{\P(V\oplus\C)}([v,1])=
\frac{\Phi_V(v)}{\pi(1+\|v\|^2)}=\const\frac{\langle\xi
v,v\rangle}{1+\|v\|^2},
\end{equation}
where $[v,1]$ denotes the line through $(v,1)$.  Denote by $\bar X$
the closure of $X$ in $\P(V\oplus\C)$.  Let $X_\infty$ be the divisor
at infinity $\bar X\cap\P V$ in $\bar X$.  The highest-weight sets of
$X$, $\bar X$ and $X_\infty$ are closely related.  If $X$ is
irreducible, then $\hw(X)$ is a submonoid of $\lplus$, so we have the
equalities
$$
\cone_{\Q}\hw(X)=\hull_{\Q}\hw(X)= \Q_{\geq0}\cdot\hw(X).
$$
(See Section \ref{subsection:cone} for the notation.)  Also,
$\Q_{\geq0}\cdot\hw(\bar X)=\cone_{\Q}\hw(\bar X)$, because $\hw(\bar
X)$ is $\Q$-convex.  As before, let $X_\infty\red$ denote the reduced
variety associated to $X_\infty$.  By Addendum
\ref{addendum:unreduced}, $\Delta(X_\infty\red)$ is a union of convex
polytopes, one for each irreducible component of $X_\infty\red$.

\begin{theorem}\label{theorem:projective}
Let $X$ be an irreducible affine $G$-variety.
\begin{enumerate}
\item\label{projectivize} The highest-weight set of $\bar X$ is the
$\Q$-convex hull of the highest-weight set of $X_\infty\red$ and the
origin in $\tplus$\upn: $\hw(\bar
X)=\hull_{\Q}\bigl(\hw(X_\infty\red)\cup\{0\}\bigr)$\upn;
\item\label{polytope} the momentum polytope of $\bar X$ is the convex
hull of the momentum set of $X_\infty\red$ and the origin in
$\tplus$\upn: $\Delta(\bar
X)=\hull\bigl(\Delta(X_\infty\red)\cup\{0\}\bigr)$\upn;
\item\label{cone} the highest-weight sets of $X$ and $\bar X$ span the
same cone\upn: $\cone_{\Q}\hw(X)=\cone_{\Q}\hw(\bar
X)$.
\end{enumerate}
\end{theorem}

\begin{proof}
\ref{projectivize}.  First I show that
\begin{equation}\label{equation:unreduced}
\hw(\bar X)= \hull_{\Q}\bigl(\hw(X_\infty)\cup\{0\}\bigr).
\end{equation}
Let $Z$ be the affine cone on $\bar X$ and let $Y$ be the affine cone
on $X_\infty$.  Let $z\colon V\times\C\to\C$ be the projection onto
the second factor and let $f$ be the restriction of $z$ to $Z$.  Then
$f$ can be regarded as an element of $\C[Z]_1$ and as such it is an
invariant, because $G$ acts trivially on the second factor.  The
coordinate ring of $Y$ is equal to $\C[Z]/(f)$.  From the
$G$-equivariant exact sequence
\begin{equation}\label{equation:exact}
\xymatrix{{(f)}\ar@{((->}[r] & {\C}[Z]\ar@{->>}[r] & {\C}[Y]}
\end{equation}
it is clear that every irreducible representation occurring in
$\C[Y]_n$ also occurs in $\C[Z]_n$.  Therefore,
$\hw(X_\infty)\subset\hw(\bar X)$.  Furthermore, $\C[Z]_1$ contains
the copy $\C[Z]_0f$ of the one-dimensional trivial representation, and
so $0\in\hw(\bar X)$.  Consequently $\hw(\bar
X)\supset\hull_{\Q}\bigl(\hw(X_\infty)\cup\{0\}\bigr)$.  Conversely,
from \eqref{equation:exact} we see that if $R_{\lambda^*}$ occurs in
$\C[Z]_n$, then it occurs in either $\C[Z]_{n-1}$ or $\C[Y]_n$.  In
other words, if $\lambda/n\in\hw(\bar X)$, then
$\lambda/(n-1)\in\hw(\bar X)$ or $\lambda/n\in\hw(X_\infty)$.  This
implies that every element of $\hw(\bar X)$ lies on the segment
joining an element of $\hw(X_\infty)$ to the origin.  Thus, $\hw(\bar
X)\subset\hull_{\Q}(\hw\bigl(X_\infty)\cup\{0\}\bigr)$.  This proves
\eqref{equation:unreduced}.

To finish the proof of \ref{projectivize} it is enough to show that
$$
\hw(X_\infty)\subset\hull_{\Q}\bigl(\hw(X_\infty\red)\cup\{0\}\bigr).
$$
I do this by showing that for every $\nu\in\hw(X_\infty)$ there exists
a rational $q\ge1$ such that $q\nu\in\hw(X_\infty\red)$.  Let $q$ be
the largest rational number such that $\mu=q\nu\in\hw(X_\infty)$.
Such a $q$ exists by Addendum \ref{addendum:sup} and is clearly
$\ge1$.  There exist $g\in\C[Z]_n^N$ and $\lambda\in\lplus$ such that
$\mu=\lambda/n$, $g$ transforms according to the weight $\lambda^*$
under the action of $T$, and $g$ is not in $(f)$.  I assert that $g$
does not vanish identically on $Y$.  For if it did, then by the
Nullstellensatz there would exist $l$ such that $g^l\in(f)$.  Write
$g^l=hf^m$ with $h\in\C[Z]$ and $m$ as large as possible; then
$h\not\in(f)$, so $h+(f)$ is a nonzero element of $\C[Y]=\C[Z]/(f)$.
Since $f$ is an invariant of degree one, the weight of $h+(f)$ is
equal to $l\lambda$ and its degree is $nl-m$.  Therefore,
$$
\frac{nlq\nu}{nl-m}=\frac{nl\mu}{nl-m}=
\frac{l\lambda}{nl-m}\in\hw(X_\infty),
$$
which contradicts the maximality of $q$.  We conclude that $g(y)\ne0$
for some $y$ in $Y$.  This means that $g$ represents a nonzero element
of $\C[Y\red]$, and hence $\mu\in\hw(X_\infty\red)$.

\ref{polytope}.  This follows immediately from \ref{projectivize} and
Theorem \ref{theorem:kaehler}.

\ref{cone}.  Note first that $\cone\hw(\bar X)=\cone\hw(Z)$, because
$Z$ is the affine cone on $\bar X$.  Furthermore, the affine
$G$-variety $Z-Y$ is saturated in $Z$, because $Y$ is defined as the
zero set of the invariant function $f$.  Therefore, $\hw(Z-Y)=\hw(Z)$
by Lemma \ref{lemma:open}.  Moreover, the $G$-equivariant map from
$V\oplus\C$ to itself sending $(x,t)$ to $(tx,t)$ maps
$X\times\C^\times$ isomorphically onto $Z-Y$.  It follows that
$\C[Z-Y]$ is isomorphic to $\C[X]\otimes\C[f,f\inv]$ as a $G$-algebra,
where $G$ acts trivially on $\C[f,f\inv]$.  This implies
$\hw(Z)=\hw(X)$.  In sum, we have shown that
$\cone\hw(X)=\cone\hw(\bar X)$.
\end{proof}

The proof shows that $\Delta(\bar X)$ is in fact equal to the
\emph{join\/} of $\Delta(X_\infty\red)$ with the origin in $\tplus$,
that is, the union of all intervals joining points in
$\Delta(X_\infty\red)$ to the origin.

Here is the main result of this section.

\begin{theorem}\label{theorem:affine}
For every $G$-stable irreducible closed affine subvariety $X$ of $V$
the set $\Delta(X)$ is equal to $\cone\hw(X)$.  In particular\upn, it
is a rational convex polyhedral cone.
\end{theorem}

\begin{proof}
Because the monoid $\hw(X)$ is finitely generated, it spans a rational
convex polyhedral cone in $\lie t^*$.  To prove that
$\Delta(X)=\cone\hw(X)$ it suffices to prove that
\begin{gather}
\label{equation:subset} \Delta(X)\subset\cone\hw(X);\\
\label{equation:quadrant}
\hull\{\lambda_1,\dots,\lambda_k\}\subset b\,\Delta(X) \quad\text{for
all $b>0$ and $\lambda_1,\dots$, $\lambda_k$ in $\hw(X)$.}
\end{gather}

\emph{Proof of\/} \eqref{equation:subset}.  By
\ref{equation:projective2}, $\Delta(X)$ is a subset of the cone on
$\Delta(\bar X)$.  By Theorem \ref{theorem:kaehler}, $\Delta(\bar X)$
is equal to the closure of $\hw(\bar X)$ in $\tplus$, so by \ref{cone}
of Theorem \ref{theorem:projective}, the cone on $\Delta(\bar X)$ is
equal to the convex hull of $\hw(X)$.  Consequently, $\Delta(X)$ is a
subset of the convex hull of $\hw(X)$.

\emph{Proof of\/} \eqref{equation:quadrant}.  Let $b$ be any positive
number, let $\mu\in\tplus$ and let $Y$ be the product $V\times K\mu$
with symplectic form $\sigma=b\omega_V+\omega_\mu$ and momentum map
$\Psi=b\Phi_V+\iota_\mu$.  I assert that
\begin{equation}\label{equation:admissible} 
\text{the momentum map $\Psi$ is admissible for all $b>0$ and
$\mu\in\tplus$.}
\end{equation}
Assuming this for the moment, let us consider the trivial line bundle
$\O_V=V\times\C$ on $V$ with the Hermitian metric defined by the
Gaussian $h(v)=\exp(-\pi b\|v\|^2)$.  The curvature form of $(\O_V,h)$
is $-2\pi\sq\,b\omega_V$.  Lift the $K$-action on $V$ to $\O_V$ by
letting $K$ act trivially on the fibre $\C$.  Then the fibre metric is
$K$-invariant, the $K$-invariant (or $G$-invariant) holomorphic
sections of $\O_V$ are just the $G$-invariant holomorphic functions on
$V$, and the associated momentum map is $b\Phi_V$.

Let us apply this to the special case $\mu=\sum_ia_i\mu_i$, where the
$a_i$ are nonnegative numbers and the $\mu_i$ are in $\hw(X)$.
Consider the varieties $Y_i=V\times K\mu^*_i$, on which we have the
line bundles $L_i=\O_V\otimes\O_{\mu^*_i}$, symplectic forms
$\sigma_i=b\omega_V+\omega_{\mu_i}$ and momentum maps
$\Psi_i=b\Phi_V+\iota_{\mu_i}$.  By \eqref{equation:admissible}, the
momentum map $\Psi=\sum_ia_ip_i^*\Psi_i$ on $Y=V\times K\mu^*$ is
admissible for all $b>0$.  Moreover, by the definition of $\hw(X)$,
\begin{align*}
\mu_i\in\hw(X) & \iff \text{there exists a $G$-equivariant linear
surjection $\C[X]\to R_{\mu_i}$} \\ & \iff \text{there exists a
nonzero $G$-invariant vector in $\C[X]\otimes R_{\mu^*_i}$} \\ & \iff
\text{there exists a nonzero $G$-invariant algebraic section of
$\O_X\otimes\O_{\mu^*_i}$}.
\end{align*}
Remark \ref{remark:saturated} now implies that for all $b>0$ the
polytope $\hull\{\mu_1,\mu_2,\dots,\mu_k\}$ is contained in
$b\Psi(X)$.

\emph{Proof of\/} \eqref{equation:admissible}.  After rescaling $\Psi$
we may assume that $b=1$.  Let $\lie F(t,\cdot)$ be the flow of
$-\grad|\Psi|^2$.  We have to show that for every $y$ in $Y$ the
trajectory $\lie F(t,y)$ is contained in a compact subset of
$Y=V\times K\mu$.  Since $K\mu$ is compact, we need only show that the
projection of $\lie F(t,y)$ onto $V$ is contained in a compact subset
of $V$.  For any function $u$ on $Y$, let $\grad_Vu$ denote the
component of $\grad u$ along $V$.  Then using
\eqref{equation:gradient} we find for every pair $(v,\beta)$ in $Y$
\begin{multline*}
\grad_V|\Psi(v,\beta)|^2=\grad_V|\Phi_V(v)+\beta|^2\\
=\grad_V\bigl(|\Phi_V(v)|^2+2(\Phi_V(v),\beta)+|\mu|^2\bigr)
=\grad_V|\Phi_V(v)|^2+2J\beta^\flat_{V,v}.
\end{multline*}
This implies
\begin{multline}\label{equation:angle}
\langle\grad_V|\Psi(v,\beta)|^2,v\rangle =4|\Phi_V(v)|^2+2\langle
v,J\beta^\flat_{V,v}\rangle =4|\Phi_V(v)|^2+2\bigl\langle
v,\grad\Phi_V^{\beta^\flat}(v)\bigr\rangle\\
=4|\Phi_V(v)|^2+4\Phi_V^{\beta^\flat}(v)
=4|\Phi_V(v)|^2+4\bigl(\Phi_V(v),\beta\bigr)
=4\biggl|\Phi_V(v)+\frac\beta{2}\biggr|^2-|\mu|^2,
\end{multline}
where I have used the fact that $\Phi^\xi$ is homogeneous of degree
two for all $\xi$ in $\lie k$, and that $|\Phi_V|^2$ is homogeneous of
degree four.  Now suppose $(v,\beta)$ is a point where
$\langle\grad_V|\Psi(v,\beta)|^2,v\rangle\leq0$.  Then it follows from
\eqref{equation:angle} that $\Phi_V(v)$ is contained in the ball of
radius $|\mu|/2$ about the point $-\beta/2\in\lie k^*$.  Therefore,
$\Phi_V(v)$ is contained in the ball of radius $|\mu|$ about the point
$-\beta$.  In other words, $\bigl|\Phi_V(v)+\beta\bigr|^2\leq|\mu|^2$,
that is, $|\Psi(v,\beta)|^2\leq|\mu|^2$.  In short,
\begin{equation}\label{equation:outwards}
\langle\grad_V|\Psi(v,\beta)|^2,v\rangle\leq0
\;\Longrightarrow\;|\Psi(v,\beta)|^2\leq|\mu|^2.
\end{equation}
Now let $\gamma(t)$ be the projection onto $V$ of the trajectory $\lie
F\bigl(t,(v,\beta)\bigr)$ through any point $(v,\beta)\in V\times
K\mu$.  It follows from \eqref{equation:outwards} that we have the
following two (non-exclusive) possibilities:
$\bigl\langle\grad_V\bigl|\Psi\bigl(\lie
F(t,(v,\beta))\bigr)\bigr|^2,v\bigr\rangle>0$ for all $t>0$, or
$|\Psi\bigl(\lie F(s,(v,\beta))\bigr)|^2\leq|\mu|^2$ for some $s>0$.
In the first case, the curve $\gamma(t)$ is trapped inside the ball of
radius $\|v\|$ about the origin in $V$.  In the second case,
$|\Psi\bigl(\lie F(t,(v,\beta))\bigr)|^2\leq|\mu|^2$ for all $t\geq
s$, because $|\Psi|^2$ is decreasing along $\lie
F\bigl(t,(v,\beta)\bigr)$.  This implies that
$|\Phi_V|^2\bigl(\gamma(t)\bigr)\leq4|\mu|^2$ for all $t\geq s$.
Moreover, $\gamma(t)$ is contained in the $G$-orbit through $v$.  It
now follows from Lemma \ref{lemma:bounded} below that $\gamma(t)$ is
contained in a compact subset of $V$.
\end{proof}

The following lemma implies that for every point $v$ in $V$ the
restriction of $\Phi_V$ to the affine variety $\overline{Gv}$ is a
proper map.

\begin{lemma}\label{lemma:bounded}
For every bounded subset $D$ of\/ $\lie k^*$ and for every bounded
subset $B$ of $V$ the intersection $\Phi_V\inv(D)\cap GB$ is a bounded
subset of $V$.
\end{lemma}

\begin{proof} 
Let $B$ be any bounded subset of $V$.  Let $\{g(n)\}_{n\geq0}$ be a
sequence of elements of $G$, let $\{v(n)\}_{n\geq0}$ be a sequence of
vectors in $B$, and put $f(n)=\|g(n)v(n)\|^2$.  Suppose that
$\lim_{n\to\infty}f(n)=\infty$.  We need to show that the sequence
$\bigl\{\Phi_V\bigl(g(n)v(n)\bigr)\bigr\}_{n\geq0}$ is
unbounded.  By the Cartan decomposition, $G=K\exp(\sq\,\lie t)K$, so
we can write $g(n)=k(n)\exp\bigl(\sq\,\xi(n)\bigr)h(n)$, where $k(n)$,
$h(n)\in K$ and $\xi(n)\in\lie t$.  Choose an orthonormal basis
$\{e_i\}$ of $V$ with respect to which the $T$-action is diagonal.
Then there exist $\beta_i\in\lie t^*$ such that $\xi
e_i=\sq\,\beta_i(\xi)e_i$ for all $\xi\in\lie t$.  Write
$h(n)v(n)=\sum_iv_i(n)e_i$ and $\rho=\sup\{\,\|v\|:v\in B\,\}$; then
\begin{gather*}
|v_i(n)|\leq\rho\quad\text{for all $i$ and $n$},\\
g(n)v=k(n)\sum_i\mathrm e^{\beta_i(\xi(n))}v_i(n)e_i,\\
f(n)=\sum_i\mathrm e^{2\beta_i(\xi(n))}|v_i(n)|^2.
\end{gather*}
Consider the set $I$ consisting of all $i$ such that the sequence
$\exp\bigl(2\beta_i(\xi(n))\bigr)|v_i(n)|^2$ is unbounded.  Then $I$ is
nonempty, because $\lim_{n\to\infty}f(n)=\infty$.  After replacing
$\{g(n)\}$ and $\{v(n)\}$ by suitable subsequences we may assume that
$$
\lim_{n\to\infty}\mathrm e^{2\beta_i(\xi(n))}|v_i(n)|^2=\infty
$$
for all $i\in I$.  Then
$\lim_{n\to\infty}\beta_i\bigl(\xi(n)\bigr)=\infty$ for all $i\in I$,
because $|v_i(n)|$ is bounded.  This implies there exists an
$\eta\in\lie t$ with $\beta_i(\eta)>0$ for all $i\in I$.  We may
assume $\eta$ has length 1.  Then $k(n)\inv\eta$ has length 1, so
\begin{equation}\label{equation:greater}
\bigl|\Phi_V\bigl(g(n)v\bigr)\bigr|^2\geq
\bigl|\Phi^{k(n)\inv\eta}\bigl(g(n)v\bigr)\bigr|^2
=\bigl|\Phi^\eta\bigl(\exp(\sq\,\xi(n))h(n)v\bigr)\bigr|^2.
\end{equation}
A straightforward computation using \eqref{equation:quadratic2} shows
that for all $\eta\in\lie t$
\begin{equation}\label{equation:equal}
\Phi^\eta\bigl(\exp(\sq\,\xi(n))h(n)v\bigr)
=\frac1{2}\sum_i\beta_i(\eta)\mathrm e^{2\beta_i(\xi(n))}|v_i(n)|^2.
\end{equation}
The vector $\eta$ was chosen in such a way that all unbounded terms in
the right-hand side of \eqref{equation:equal} tend to $\infty$ for
$n\to\infty$.  It now follows from \eqref{equation:greater} that
$\bigl|\Phi_V\bigl(g(n)v\bigr)\bigr|^2$ tends to $\infty$ for
$n\to\infty$.
\end{proof}

\begin{corollary}\label{corollary:infinity}
The momentum cone of $X$ is the cone over the momentum polytope of the
projectivization of $X$\upn: $\Delta(X)=\cone\Delta(\bar
X)=\Q_{\geq0}\cdot\Delta(\bar X)$.
\end{corollary}
 
\begin{proof}
Combine Theorem \ref{theorem:projective}.\ref{cone} with Theorem
\ref{theorem:affine}.
\end{proof}

\begin{corollary}\label{corollary:embedding}
The momentum cone of $X$ is a closed subset of $\tplus$\upn, and it
does not depend on the embedding of $X$ into the unitary $K$-module
$V$.
\qed
\end{corollary}

\begin{corollary}\label{corollary:fibre}
If $X$ is normal\upn, all fibres of the momentum map $\Phi_X$ are
connected.
\end{corollary}

\begin{proof}
The function $|\Psi|^2=|\Phi_V+\iota_\mu|^2$ is real-algebraic on
$Y=V\times K\mu$ and has therefore only finitely many critical levels.
This implies that the Morse decomposition of $Y$ with respect to
$|\Psi|^2$ is finite.  By the proof of Theorem \ref{theorem:affine},
the momentum map $\Psi$ is admissible.  It now follows from the
results quoted in Section \ref{subsection:flow} that its zero level,
which is $K\Phi_V\inv(-\mu)$, is a deformation retract of an open
subset of $V$ the complement of which is a finite union of
complex-analytic subsets of positive codimension.  The same holds with
$V$ replaced by $X$, because the flow of $|\Psi|^2$ leaves $X\times
K\mu$ invariant.  Since $X$ is normal, the complement of a finite
number of analytic subsets is always connected.  This implies that
$K\Phi_X\inv(-\mu)$, and hence $\Phi_X\inv(-\mu)$, are connected.
\end{proof}

\begin{remark}
The zero fibre $\Phi_X\inv(0)$ is connected regardless of whether $X$
is normal.  The reason is that the function $|\Phi_V|^2$ has only one
critical level, namely $0$.
\end{remark}

\begin{corollary}\label{corollary:saturated}
Let $Y$ be a saturated Zariski-open subvariety of $X$.  Then
$\Delta(Y)=\Delta(X)$.
\end{corollary}

\begin{proof}
Evidently, $\Delta(Y)$ is contained in $\Delta(X)$.  For the reverse
inclusion it suffices to show that
$\hull\{\lambda_1,\dots,\lambda_k\}\subset b\,\Delta(Y)$ for all $b>0$
and $\lambda_1,\dots$, $\lambda_k$ in $\hw(X)$.  The proof of this
fact is a straightforward generalization of the proof of
\eqref{equation:quadrant}.  (Cf.\ Remark \ref{remark:saturated}.)
\end{proof}

\begin{example}\label{example:real}
Suppose the affine $G$-variety $X$ is defined over the real numbers in
the sense that the complex $G$-algebra $\C[X]$ is the complexification
of a real $K$-algebra of finite type.  Then complex conjugation
defines an antilinear involution on the $G$-module $\C[X]$, so
whenever an irreducible representation $R_\lambda$ occurs in $\C[X]$,
its contragredient representation $R_\lambda^*$ also occurs.  It
follows from this that the monoid $\hw(X)$ is invariant under the
involution $*\colon\lplus\to\lplus$.  Therefore, by Theorem
\ref{theorem:affine}, the cone $\Delta(X)$ is invariant under the
involution $*\colon\tplus\to\tplus$.  This can also be shown directly
as follows.  We may assume the embedding of $X$ into the $G$-module
$V$ to be defined over the real numbers (in the sense that both $V$
and the $G$-morphism $X\to V$ are defined over the reals).  Then $X$
is invariant under complex conjugation on $V$.  From
\eqref{equation:quadratic1} one deduces immediately that
$\Phi_V^\xi(\bar v)=-\Phi_V^\xi(v)$.  Hence $\Phi_V(X)=-\Phi_V(X)$,
and therefore $\Delta(X)^*=-w_0\Delta(X)=\Delta(X)$.
\end{example}

\begin{example}[Peter-Weyl]\label{example:peterweyl}
The group $G$ is an affine variety in its own right, and it acts on
itself by left multiplication: $\lie L_gh=gh$, and by right
multiplication: $\lie R_gh=hg\inv$.  Consider the $\lie L\times\lie
R$-action of $G\times G$ on $G$.  Let us denote the highest-weight set
of $G$ for this action by $\hw(G,\lie L\times\lie R)$ and the momentum
cone (with respect to any algebraic $G\times G$-equivariant embedding
of $G$ into a unitary $K\times K$-module) by $\Delta(G,\lie
L\times\lie R)$.  The monoid of dominant weights of $G\times G$ is
simply the product $\lplus\times\lplus$ and its positive Weyl chamber
is $\tplus\times\tplus$.  By the Peter-Weyl Theorem the coordinate
ring of $G$ is a direct sum of irreducible $G\times G$-modules:
$\C[G]=\bigoplus_{\lambda\in\lplus} R_\lambda\otimes R^*_\lambda$.
This implies that the highest-weight monoid $\hw(G,\lie L\times\lie
R)$ is equal to the subset
$\{\,(\lambda,\lambda^*):\lambda\in\lplus\,\}$ of
$\lplus\times\lplus$.  By Theorem \ref{theorem:affine}, the momentum
cone $\Delta(G,\lie L\times\lie R)$ is therefore the ``anti-diagonal''
$\{\,(\mu,\mu^*):\mu\in\tplus\,\}$ inside $\tplus\times\tplus$.
Notice that this set is $*$-invariant as it should be, because $G$ is
defined over the real numbers.
\end{example}

\begin{example}\label{example:leftright}
We use the notation of the previous example.  There are three different
embeddings of $G$ into $G\times G$: the maps $i_1(g)=(g,1)$,
$i_2(g)=(1,g)$ and $d(g)=(g,g)$.  Pulling back the $\lie L\times\lie
R$-action via these three embeddings yields three actions of $G$ on
itself: the actions $\lie L$ (left multiplication), $\lie R$ (right
multiplication) and $\lie C$ (conjugation).  The momentum maps for
these actions are obtained by composing the $\lie L\times\lie
R$-momentum map with the maps $i_1^*$, $i_2^*$ and $d^*$,
respectively.  Thus we find that $\Delta(G,\lie L)$ and
$\Delta(G,\lie R)$ are equal to the positive Weyl chamber, $\tplus$,
and $\Delta(G,\lie C)$ is the positive Weyl chamber of the semisimple
part of $K$: $\Delta(G,\lie C)=\tplus\cap[\lie k,\lie k]$.
\end{example}

\begin{example}[Gelfand's variety $G\qu N$]\label{example:gelfand}
The action $\lie L$ of $G$ on itself descends to an action of $G$ on
$G\qu N$.  Every irreducible $G$-module occurs exactly once in the
coordinate ring $R=\C[G]^N$, so, once again, $\hw(G\qu N)=\lplus$ and
$\Delta(G\qu N)=\tplus$.

We can embed $G\qu N$ into affine space and compute the momentum
polytope of its projectivization.  First assume $G$ is semisimple and
simply connected.  Then the algebra $R$ is generated by the subspace
$E=\bigoplus_{i=1}^rR_{\pi_i}$, where $\pi_1$, $\pi_2,\dots$, $\pi_r$
are the fundamental weights of $G$ and $r$ is the rank of $G$.  Choose
a highest-weight vector $v_i$ in each of the $R_{\pi_i}$.  Consider
the left-$G$-equivariant map from $G$ to $E$ defined by sending the
identity of $G$ to the vector $v_1\oplus v_2\oplus\cdots\oplus v_r$.
This map is right-$N$-equivariant, so it descends to a map from $G\qu
N$ to $E$, which is by construction an embedding.  Let us identify
$G\qu N$ with its image in $E$.  It is not hard to show that the
subvariety $G\qu N$ is invariant under the standard $\C^\times$-action
on $E$, and the divisor at infinity $(G\qu N)_\infty$ is therefore the
quotient $(G\qu N-\{0\})/\C^\times$.  In other words, $G\qu N$ is the
affine cone on $(G\qu N)_\infty$.  It now follows immediately from
Theorem \ref{theorem:kaehler} that the momentum polytope of $(G\qu
N)_\infty$ (with respect to any $K$-invariant inner product on $E$) is
the $r-1$-dimensional simplex spanned by the fundamental weights.  By
Theorem \ref{theorem:projective}, the momentum polytope of the
projective closure of $G\qu N$ is therefore the $r$-dimensional
simplex spanned by the fundamental weights and the origin in $\tplus$.

Now assume $G$ is a torus of dimension $k$.  Then the subgroup $N$ is
trivial and so $R=\C[G]$ and $G\qu N=G$.  Let $\zeta_1$,
$\zeta_2,\dots$, $\zeta_k$ be a basis over $\Z$ of the weight lattice
$\Lambda^*$, and identify $\lie t^*$ with $\R^k$ by sending this basis
to the standard basis in $\R^k$.  This choice of basis gives an
identification of $G$ with the product $(\C^\times)^k=
\{\,(t_1,t_2,\dots,t_k):t_i\in\C^\times\,\}$.  A \emph{closed\/}
affine embedding of $G$ is given by sending $(t_1,t_2,\dots,t_k)$ to
$(t_1,t_1\inv,t_2,t_2\inv, \dots,t_k,t_k\inv)\in\C^{2k}$.  The
projective closure of $G$ in $\P^{2k}$ is a product of $k$ copies of
$\P^1$.  The divisor at infinity $G_\infty$ contains $2^k$ fixed
points for the action of $G$, whose images under the momentum map are
the points $\pm\zeta_1\pm\zeta_2\pm\dots\pm\zeta_k$.  Theorem
\ref{theorem:projective} now implies that the momentum polytope of the
projective closure of $G$ is the parallelepiped spanned by these $2^k$
points.

For an arbitrary connected reductive group $G$, the variety $G\qu N$
can be embedded into affine space in a similar way, by choosing a
basis of the monoid of highest weights.  One can show that the
momentum polytope of the projective closure of $G\qu N$ under such an
embedding is the product of the simplex spanned by the origin and the
fundamental weights of $[\lie k,\lie k]$, and the parallelepiped
spanned by the points $\pm\zeta_1\pm\zeta_2\pm\dots$, where the
$\zeta_i$ are a basis of the weight lattice of $\lie z(\lie k)$, the
centre of $\lie k$.
\end{example}

\begin{example}[associated bundles]\label{example:associated}
Let $F$ be a reductive subgroup of $G$ and let $Y$ be an affine
$F$-variety.  Consider the bundle $X=G\times^FY$ associated to the
principal fibration $F\to G\to G/F$.  The action $\lie L$ of $G$ on
itself induces a $G$-action on $X$.  Also, $X$ is an affine variety
with coordinate ring $\C[X]=(\C[G]\otimes\C[Y])^F$.  Note that if
$F^0$ is the identity component of $F$, there is a finite map
$G\times^{F^0}Y\to X$, so $\Delta(X) =\Delta(G\times^{F^0}Y)$ by Lemma
\ref{lemma:finite}.  This means we may assume $F$ to be connected.
The categorical quotient of $X$ by $N$ has coordinate ring
$\C[X]^N=(\C[G]\otimes\C[Y])^{F\times N}$, which is isomorphic to
$(R\otimes\C[Y])^F$, where $R$ is the ring $\C[G]^N$, on which $F$
acts by \emph{left\/} multiplication.  By Lemma \ref{lemma:unipotent}
and Theorem \ref{theorem:affine}, the momentum cone of $X$ is
therefore the convex cone spanned by the weights of the action of the
maximal torus $H$ on the algebra $(R\otimes\C[Y])^F$ defined by right
multiplication on $R$.  In general, this is hard to calculate
explicitly.
\end{example}

\begin{example}[tori]\label{example:torus}
In the setting of the previous example, let us assume that $G=H$ is a
torus, and let us write $F=H_1$.  As noted above, we may assume $H_1$
to be connected.  Then $H_1$ is the complexification of a subtorus
$T_1$ of $T$.  We let $T_2$ be the quotient $T/T_1$ and identify it
with a complement of $T_1$ in $T$, so that $T\cong T_1\times T_2$.  Put
$H_2=(T_2)\co$.  Then $H\cong H_1\times H_2$ and
$$
X=H\times^{H_1}Y\cong(H_1\times H_2)\times^{H_1}Y =H_2\times
(H_1\times^{H_1}Y)=H_2\times Y.
$$
Therefore, by Remark \ref{remark:product}, 
$$
\hw(X)\cong\hw(H_2)\times\hw(Y)=\Lambda^*_2\times\hw(Y),
$$
where $\Lambda^*_2$ is the weight lattice of $H_2$.  Consequently,
$\Delta(X)\cong\lie t_2^*\times\Delta(Y)$.  (These identifications
depend on the splitting $H\cong H_1\times H_2$.  An invariant way of
stating these facts is: $X$ is a trivial principal $H_2$-bundle over
$Y$; $\hw(X)$ is equal to the preimage of $\hw(Y)$ under the canonical
projection $\Lambda^*\to\Lambda^*_1$; and $\Delta(X)$ is equal to the
preimage of $\Delta(Y)$ under the canonical projection $\lie
t^*\to\lie t_1^*$.)  If $Y$ is a vector space, then by
\eqref{equation:weights}, $\Delta(X)\cong\lie
t_2^*\times-\cone\{\nu_1,\dots,\nu_l\}$, where $\nu_1,\dots$,
$\nu_l$ are the weights of the $H_1$-action on $Y$.
\end{example}

\subsection{The momentum cone and \'etale slices}
\label{subsection:etale}

Remarkably, the momentum cone $\Delta(X)$ of an affine $G$-variety $X$
turns out to be entirely determined by infinitesimal data at any point
on a closed $G$-orbit.  I shall deduce this from Luna's \'etale slice
theorem.  First I discuss a variation on Lemma \ref{lemma:finite}.

\begin{proposition}\label{proposition:finite}
Let $X$ and $Y$ be affine $G$-varieties\upn, let $\phi\colon X\to Y$
be a $G$-morphism\upn, let $x$ be a point in $X$\upn, and let
$y=\phi(x)$.  Suppose that $\phi$ has finite fibres\upn, that the image
of $\phi$ is open in $Y$\upn, and that the orbits $Gx$ and $Gy$ are
closed.  Then $\hw(Y)$ is contained in $\hw(X)$\upn, and $\hw(X)$ is
contained in the cone on $\hw(Y)$.
\end{proposition}

\begin{proof}
The orbit $Gy$ and the complement of $\phi(X)$ are $G$-stable
Zariski-closed subsets of $Y$.  Because $G$-invariant polynomial
functions separate $G$-stable Zariski-closed subsets, there exists an
$f\in\C[Y]^G$ that vanishes outside $\phi(X)$ and satisfies $f(y)=1$.
Let $Y'=Y_f$ and $X'=\phi\inv(Y)$.  Then $X'$ and $Y'$ are saturated
affine open subsets of $X$, resp.\ $Y$, containing the orbits $Gx$,
resp.\ $Gy$, and the restriction of $\phi$ to $X'$ is surjective onto
$Y'$.  Then $\hw(X')=\hw(X)$ and $\hw(Y')=\hw(Y)$ by Lemma
\ref{lemma:open}.  We are therefore reduced to proving that
$\hw(Y')\subset\hw(X')$ and $\hw(X')\subset\cone\hw(Y')$.

Let $\lie C$ be the integral closure of $\C[Y']$ in $\C[X']$, and let
$Z=\spec\lie C$.  By Luna's equivariant version of Zariski's Main
Theorem (\cite{lu:sl}, part~I), the natural maps $\iota\colon X'\to Z$
and $\psi\colon Z\to Y'$ have the following properties: $\iota$ is an
open immersion, $\psi$ is a finite morphism, and
$\phi=\psi\circ\iota$.  Also, $\psi$ is surjective, because
$\phi\colon X'\to Y'$ is.  Consequently, $\hw(Y')\subset\hw(Z)$ and
$\hw(Z)\subset\cone\hw(Y')$ by Lemma \ref{lemma:finite}.  So if we can
show that $\hw(X')=\hw(Z)$, we are done.  Let us identify $X'$ with
its image $\iota(X')$ in $Z$.  The orbit $Gx\subset X'$ is closed in
$Z$ (cf.\ \cite{lu:sl}, p.~94): since $Gy$ is closed in $Y'$ and
$\psi$ is finite, $\psi\inv(Gy)$ is closed in $Z$ and consists of a
finite number of orbits, one of which is $Gx$.  The conclusion is that
$Gx$ and the complement of $X'$ in $Z$ are $G$-stable closed subsets
of $Z$.  This implies the existence of a $G$-invariant $h\in\C[Z]$
that vanishes outside $X'$ and satisfies $h(x)=1$.  Then $X'_h$ is a
$G$-stable affine open subset of $X'$, and it is saturated as a subset
of both $X'$ and $Z$.  Hence, by Lemma \ref{lemma:open},
$\hw(X')=\hw(X'_h)=\hw(Z)$.
\end{proof}

\begin{theorem}\label{theorem:level}
Let $X$ be an affine $G$-variety\upn, let $x$ be a point on a closed
$G$-orbit\upn, and let $S_x$ be an \'etale slice at $x$.  Then the
momentum cone of $X$ is equal to the momentum cone of
$G\times^{G_x}S_x$.  If $x$ is a smooth point of $X$\upn, then
$\Delta(X)=\Delta\bigl(G\times^{G_x}V_x\bigr)$\upn, where $V_x$ is the
tangent space to $S_x$ at $x$.
\end{theorem}

\begin{proof}
By Luna's Etale Slice Theorem the natural map from the bundle
$G\times^{G_x}S_x$ into $X$ is \'etale and its image is Zariski-open.
Furthermore the $G$-orbits through the point $[1,x]$ in
$G\times^{G_x}S_x$ and the point $x$ in $X$ are closed.  It now follows
from Proposition \ref{proposition:finite} that $\hw(X)$ and
$\hw\bigl(G\times^{G_x}S_x\bigr)$ span the same cone.  Hence
$\Delta(X)=\Delta\bigl(G\times^{G_x}S_x\bigr)$ by Theorem
\ref{theorem:affine}.

If $x$ is a smooth point of $X$, we may assume the \'etale slice $S_x$
to be smooth, and there exists a $G_x$-equivariant \'etale morphism
$\psi\colon S_x\to V_x$ with Zariski-open image.  The map $\psi$
extends to a $G$-equivariant map $G\times^{G_x}S_x\to
G\times^{G_x}V_x$, which is \'etale and has Zariski-open image as
well.  Again by Proposition \ref{proposition:finite},
$\hw\bigl(G\times^{G_x}S_x\bigr)$ and
$\hw\bigl(G\times^{G_x}S_x\bigr)$ span the same cone.  We conclude that
$\Delta(X)=\Delta\bigl(G\times^{G_x}V_x\bigr)$.
\end{proof}

\begin{corollary}\label{corollary:samecone}
The cone on $\hw\bigl(G\times^{G_x}S_x\bigr)$ is independent of the
point $x$.  Here $x$ ranges over the set of all points in $X$ through
which the $G$-orbit is closed.
\qed
\end{corollary}

The following result is a necessary condition for the origin to be an
extreme point of $\Delta(X)$.  Here $[G,G]$ denotes the commutator
subgroup of $G$.  Note that for every subgroup $F$ of $G$, $[G,G]F$ is
a subgroup of $G$, because $[G,G]$ is normal.  If $F$ is a closed
reductive subgroup, then so is $[G,G]F$.

\begin{theorem}\label{theorem:propercone}
Assume that $\Delta(X)$ is a proper cone.  Then for every point $x$
such that $Gx$ is closed the following condition holds\upn:
$G=[G,G]G_x$.
\end{theorem}

\begin{proof}
Let $x$ be any point such that $Gx$ is closed.  Let $Y$ denote the
homogeneous space $G/G_x$ and $Z$ the homogeneous space $G/[G,G]G_x$.
Consider the maps
$$
\xymatrix{X & Y\ar[r]^-\tau\ar[l]_-\iota & Z,}
$$
where $\iota$ is the $G$-map sending the coset $1G_x$ to $x$ and
$\tau$ is the canonical projection.  Clearly, $\Delta(Y)$ is a subset
of $\Delta(X)$ and, by Remark \ref{remark:surjective} and Theorem
\ref{theorem:affine}, $\Delta(Z)$ is a subset of $\Delta(Y)$.
Therefore, since $\Delta(X)$ is a proper cone, so is $\Delta(Z)$.  On
the other hand, the torus $G/[G,G]$ acts transitively on $Z$, so
$\Delta(Z)$ is a vector space.  (Cf.\ Example \ref{example:torus}.)
We conclude that $Z$ is a point, in other words, $G=[G,G]G_x$.
\end{proof}

Note that if $G$ is semisimple, the condition $G=[G,G]G_x$ is void.
This is as it should be, because in this case the positive Weyl
chamber $\tplus$ is a proper cone, so every cone contained in it is a
proper cone.

\section{Stein varieties}\label{section:stein}

In this section I prove a convexity theorem for certain Stein
$K$-varieties, Theorem \ref{theorem:stein}.  It can be regarded as a
local version of Theorem \ref{theorem:affine}.  The results are far
from optimal, but will be sufficient for our purposes.  Let me start
with a number of elementary observations on K\"ahler potentials and
momentum maps.

\begin{lemma}\label{lemma:potential}
Suppose $Y$ is a connected complex manifold and $\rho$ a strictly
plurisubharmonic function on $Y$.  Let $\sigma$ be the K\"ahler form
$\sq\,\partial\bar\partial\rho$ with associated Riemannian metric
$\langle\cdot,\cdot\rangle$\upn, and let $\vartheta$ be the
Hamiltonian vector field of $\rho$.  Then the vector field
$J\vartheta=\grad\rho$ is expanding\upn:
$\Lie_{J\vartheta}\sigma=2\sigma$.
\end{lemma}

\begin{proof}
Let $J\colon TY\to TY$ denote the complex structure on $Y$ and also
the transpose operator $T^*Y\to T^*Y$.  For all functions $f$,
\begin{equation}\label{equation:ddbar}
\bar\partial f=\frac1{2}(d+\sq\,Jd)f\quad\text{and}\quad
dJdf=-2\sq\,d\bar\partial f=-2\sq\,\partial\bar\partial f.
\end{equation}
Moreover, for all functions $f$ and all tangent vectors $\eta$, 
$$
\sigma(\grad f,\eta)=-\langle\grad
f,J\eta\rangle=-df(J\eta)=-Jdf(\eta),
$$
and therefore $Jdf=-\iota_{\grad f}\sigma$.  Together with
\eqref{equation:ddbar} this implies that $\Lie_{\grad\rho}\sigma
=d\iota_{\grad\rho}\sigma=-dJd\rho
=2\sq\,\partial\bar\partial\rho=2\sigma$.
\end{proof}

Let me add to this that the vector fields $\vartheta$ and $J\vartheta$
are usually not holomorphic.

The symplectic form on $Y$ being exact, every symplectic action of $K$
on $Y$ has a momentum map $\Psi$.  It turns out that the flow of the
vector field $-\grad\rho$ has the peculiar property that it pushes
forward under $\Psi$ to a flow on $\lie k^*$, which retracts the image
of $\Psi$ exponentially to a single point in $\lie z(\lie k)^*$, where
$\lie z(\lie k)$ is the centre of $\lie k$.

\begin{proposition}\label{proposition:hamilton}
Let $Y$ and $\rho$ be as in Lemma \ref{lemma:potential}.  Suppose that
$K$ acts holomorphically on $Y$\upn, leaving $\rho$ invariant.  Put
$\alpha=-\sq\,\bar\partial\rho$ and $\Psi^\xi=\iota_{\xi_Y}\alpha$ for
all $\xi\in\lie k$.  Let $\lie G(t,\cdot)=\lie G_t(\cdot)$ denote the
flow of $-J\vartheta=-\grad\rho$.  Then

\begin{enumerate}
\item\label{primitive} the functions $\Psi^\xi$ are the components of
an equivariant momentum map $\Psi$ for the $K$-action on $Y$ with
respect to the symplectic form $\sigma$\upn;
\item\label{character} there exists a character $c$ of\/ $\lie k$ such
that $\Lie_{J\vartheta}\Psi=2(\Psi+c)$ \upn(cf.\ \cite{he:ka}\upn,
\S~$3$\upn).  Therefore
\begin{equation}\label{equation:flow}
\lie G_t^*\Psi=\mathrm e^{-2t}\Psi+(\mathrm e^{-2t}-1)c
\end{equation}
for all $t$ such that $\lie G_t$ is defined.  It follows that $\lie
G_t$ maps fibres of $\Psi$ to fibres of $\Psi$.  Moreover\upn, if
$\lie G_t$ is defined for all $t\geq0$\upn, then
$\lim_{t\to\infty}\lie G_t^*\Psi$ is the constant map $-c$\upn;
\item\label{limit} if $c=0$\upn, then the critical set of $\rho$ is
contained in the fibre $\Psi\inv(0)$.  The converse holds if $\rho$
has at least one critical point.
\end{enumerate}
\end{proposition}

\begin{proof}
Note first that $\sigma=-d\alpha$.  Further, since $\rho$ is
$K$-invariant and $K$ acts holomorphically, $\alpha$ is $K$-invariant.
Consequently,
$$
d\Psi^\xi=
d\iota_{\xi_Y}\alpha=\Lie_{\xi_Y}\alpha-\iota_{\xi_Y}d\alpha
=-\iota_{\xi_Y}d\alpha=\iota_{\xi_Y}\sigma,
$$
so $\xi_Y$ is the Hamiltonian vector field of the function $\Psi^\xi$.
An easy calculation shows that
$\{\Psi^\xi,\Psi^\eta\}=\Psi^{[\xi,\eta]}$, so $\Psi$ is
$K$-equivariant.  This proves \ref{primitive}.

Note that since $\vartheta$ is the Hamiltonian vector field of the
$K$-invariant function $\rho$, the induced vector field $\xi_Y$
commutes with $\vartheta$ for all $\xi\in\lie k$.  Being holomorphic,
$\xi_Y$ therefore commutes with $J\vartheta$ as well, so by
\ref{primitive}
$$
d\Lie_{J\vartheta}\Psi^\xi =\Lie_{J\vartheta}d\Psi^\xi
=\Lie_{J\vartheta}\iota_{\xi_Y}\sigma
=(\iota_{[J\vartheta,{\xi_Y}]}+\iota_{\xi_Y}\Lie_{J\vartheta})\sigma
=\iota_{\xi_Y}\Lie_{J\vartheta}\sigma =\iota_{\xi_Y}2\sigma
=2d\Psi^\xi.
$$
This implies the function $\Lie_{J\vartheta}\Psi^\xi-2\Psi^\xi$ is a
constant, say $2c(\xi)$, for all $\xi\in\lie k$.  It is evidently
linear in $\xi$.  From the equivariance of $\Psi$ and the fact that
$[J\vartheta,{\xi_Y}]=0$ it is now easy to deduce that
$c([\xi,\eta])=0$ for all $\xi$ and $\eta$ in $\lie k$.  This proves
the first assertion in \ref{character}.  Integrating the equation
$\Lie_{J\vartheta}\Psi=2(\Psi+c)$ yields \eqref{equation:flow}.  The
last two assertions are obvious.

Assume $c=0$.  Let $y$ be a critical point of $\rho$.  Then
\ref{character} implies that $2\Psi(y)=\bigl(\Lie_{\grad
\rho(y)}\Psi\bigr)(y)=(\Lie_0\Psi)(y)=0$, so $y\in\Psi\inv(0)$.
Conversely, assume the critical set of $\rho$ is nonempty and is
contained in the fibre $\Psi\inv(0)$.  Let $y$ be a critical point of
$\rho$.  Then from \eqref{equation:flow} we obtain
$$
0=\Psi(y)=\Psi\bigl(\lie G(t,y)\bigr)=\mathrm e^{-2t}\Psi(y)+(\mathrm
e^{-2t}-1)c=(\mathrm e^{-2t}-1)c 
$$
for $t\geq0$, and so $c=0$.
\end{proof}

The set-up of this proposition is functorial in the following
sense.  Let $Z$ be a $K$-invariant closed complex submanifold of $Y$
and let $\rho_Z=\rho|_Z$ be the restriction of $\rho$ to $Z$.  Put
$\alpha_Z= -\sq\,\bar\partial\rho$, $\sigma_Z=-d\alpha$ and
$\Psi^\xi_Z= \iota_{\xi_Z}\alpha_Z$ for $\xi\in\lie k$.  Then
$\alpha_Z=\alpha|_Z$, $\sigma_Z=\sigma|_Z$ and $\Psi_Z=\Psi|_Z$.  Of
course, the Hamiltonian vector field $\vartheta_Z$ of $\rho_Z$ is not
the restriction of $\vartheta$ to $Z$, unless $\vartheta$ happens to
be tangent to $Z$.

These observations apply to the pair of manifolds $Y=V$ and $Z=X$,
where $V$ is a $G$-representation space with a $K$-invariant inner
product as in Section \ref{section:affine}, and $X$ a $G$-stable
closed nonsingular algebraic subvariety of $V$.  We take $\rho$ to be
the function $\rho(v)=\|v\|^2/2$.  Clearly, $\sigma=
\sq\,\partial\bar\partial\rho$ is the standard symplectic form
$\omega_V$, $\Psi$ is the quadratic momentum map $\Phi_V$ given by
\eqref{equation:quadratic2}, and $J\vartheta=\grad\rho$ is the radial
vector field $v\,\partial/\partial v$ on $V$.  The idea to use the
length function as a tool in invariant theory is due to Kempf and Ness
\cite{ke:le}.  I shall frequently refer to their main result (see also
\cite{sc:to}):

\begin{theorem}\label{theorem:kempfness}
For all $v$ in $V$ the following conditions are equivalent\upn:
\begin{enumerate}
\item the orbit $Gv$ is closed\upn;
\item the restriction of $\rho$ to $Gv$ has a stationary point\upn;
\item $Gv$ intersects the zero level set of the momentum map $\Phi_V$.
\end{enumerate}
If $v$ is a stationary point of $\rho|_{Gv}$\upn, then\upn:
$\rho|_{Gv}$ takes on its minimum at $v$\upn; for all $w\in Gv$\upn,
$\rho(w)=\rho(v)$ implies $w\in Kv$\upn; and $G_v=(K_v)\co$.
\qed
\end{theorem}

Note that $\grad\rho$ is tangent to the subvariety $X$ only if $X$ is
invariant under the standard $\C^\times$-action on $V$.  Because
$X$ is closed, the restriction of $\rho$ to $X$, $\rho_X$, is a proper
function, so the forward trajectories of $-J\vartheta_X=-\grad\rho_X$
are bounded and the flow $\lie G_X(t,\cdot)$ is defined for all
$t\geq0$.   Since $\rho_X$ is real-analytic, $\lim_{t\to\infty}\lie
G_X(t,x)$ exists for all $x$ in $X$.  The properness of $\rho_X$
implies that the the flow retracts the stable set of every critical
level continuously onto the critical set.  (See Section
\ref{subsection:flow}.)   Furthermore, $\rho_X$ always has critical
points, for example minima.   If $x$ is a critical point of $\rho_X$,
it is a critical point of the restriction of $\rho_X$ to the orbit
$Gx$, and therefore $\Phi_X(x)=0$ by Theorem \ref{theorem:kempfness}.
Hence, the character $c$ in Proposition \ref{proposition:hamilton} is
$0$, so that
\begin{equation}\label{equation:lengthflow}
\Lie_{J\vartheta_X}\Phi_X= 2\Phi_X\qquad\text{and}\qquad(\lie
G_X)_t^*\Phi_X=\mathrm e^{-2t}\Phi_X.
\end{equation}

\begin{theorem}\label{theorem:stein}
Let $X$ be a $G$-stable closed nonsingular algebraic subvariety of
$V$.  Suppose that $\rho_X$ has a unique critical level.  Let $\mathcal
U$ be a basis of neighbourhoods \upn(in the classical topology on
$X$\upn) of the critical set of $\rho_X$.  Then the sets
$\Delta(U)$\upn, where $U\in\mathcal U$\upn, form a basis of
neighbourhoods of the vertex\/ $0$ of the momentum cone $\Delta(X)$.
In particular\upn, the cone spanned by $\Delta(U)$ is equal to
$\Delta(X)$ for every $U\in\mathcal U$.
\end{theorem}

\begin{proof}
Let $B_\eps$ denote the closed ball of radius $\eps$ about the origin
in $V$ and let $\delta =\min\{\,\|x\|:x\in X\,\}$ be the distance from
$X$ to the origin.  The assumption on $\rho_X$ implies that its only
critical level is the global minimum, $\delta^2/2$.  Therefore, the
sets $X\cap B_\eps$, where $\eps>\delta$, are a basis of
neighbourhoods of the critical set $X\cap B_\delta$ of $\rho_X$.  So it
suffices to prove that the sets $\Delta(X\cap B_\eps)$, where
$\eps>\delta$, form a basis of neighbourhoods of the vertex $0$ of the
momentum cone $\Delta(X)$.  The proof has three parts: first I show
that for some $\eps>\delta$ the set $\Delta(X\cap B_\eps)$ is a
neighbourhood of the vertex in $\Delta(X)$.  Then I show that the same
is true for every $\eps>\delta$.  Lastly, I prove that for every ball
$D$ about the origin in $\lie k^*$ there exists an $\eps>\delta$ such
that $\Delta(X\cap B_\eps)\subset D$.

\emph{Part\/} 1.  The flow $\lie G_X$ extends to a deformation
retraction
\begin{equation}\label{equation:deformation}
\bar{\lie G}_X\colon X\times[0,\infty]\longrightarrow X
\end{equation}
of $X$ onto the set $X\cap B_\delta\subset\Phi_X\inv(0)$.  Now take
any $\eta>\delta$.  Then for every $x$ in $X$ the trajectory $\lie
G_X(t,x)$ is contained in $B_\eta$ for sufficiently large $t$.
Moreover, by \eqref{equation:lengthflow}, $\Phi_X\bigl(\lie
G_X(t,x)\bigr)=\mathrm e^{-2t}\Phi_X(x)$.  This implies that
\begin{equation}\label{ray} 
\text{the cone spanned by $\Delta(X\cap B_\eta)$ is the whole of
$\Delta(X)$.} 
\end{equation}
Now consider the subset $S=G\cdot(X\cap B_\eta)$ of $X$.  I assert that
\begin{equation}\label{closure} 
\text{for every $x\in S$ the affine variety $\overline{Gx}$ is
contained in $S$.}
\end{equation}
Indeed, take any $x$ in $S$ and any $y$ in $\overline{Gx}$.  We have
to show that $y$ is in $S$.  Let $\lie F_t(\cdot)$ be the gradient
flow of the function $-|\Phi_V|^2$.  The limit map $\lie
F_\infty=\lim_{t\to\infty}$ is continuous, it retracts $S$ onto
$\Phi_X\inv(0)\cap B_\eta$, and it retracts $\overline{Gx}$ onto the
$K$-orbit $K(\lie F_\infty x)\subset B_\eta$.  Moreover, $\rho_X$ is
decreasing along the flow lines, so the restriction of $\rho_X$ to
$\overline{Gx}$ takes on its minimum at $\lie F_\infty(x)$.  (See
\cite{ne:to} and \cite{sc:to}.)  This implies $\lie F_\infty(y)$ is in
the $K$-orbit through $\lie F_\infty(x)$.  There are two
possibilities: either $\lie F_\infty(x)$ is in the interior of the
ball $B_\eta$, or it is on the boundary.  In the first case, the
$G$-orbit $Gy$ intersects the interior of $B_\eta$, so $y\in
G\cdot(X\cap B_\eta)=S$.  In the second case, since $\lie F_\infty(x)$
is the point closest to the origin on $\overline{Gx}$ and by
assumption the orbit $Gx$ intersects $B_\eta$, we see that $\lie
F_\infty(x)$ lies on $Gx$.  By Theorem \ref{theorem:kempfness} every
$G$-orbit intersecting $\Phi_V\inv(0)$ is closed, and therefore
$y\in\overline{Gx}=Gx\subset S$.  This proves \ref{closure}.

Next, I assert that
\begin{equation}\label{smallcone} 
\Delta(S)=\Delta(X).
\end{equation}
To see this, let $\lambda$ be any point in $\Delta(X)$.  Then
$b\lambda\in\Delta(X\cap B_\eta)$ for some $b>0$ by \ref{ray}.  Take
$x\in X\cap B_\eta$ such that $\Phi_X(x)=b\lambda$.  Then
$\Delta\bigl(\overline{Gx}\bigr)$ contains the ray through $b\lambda$
by Theorem \ref{theorem:affine}.  But $\overline{Gx}\subset S$ by
\ref{closure}, so
$\lambda\in\Delta\bigl(\overline{Gx}\bigr)\subset\Delta(S)$.  This
proves \ref{smallcone}.

Now let $D$ be any ball about the origin in $\lie k^*$.  Then
$\Phi_V\inv(D)\cap S$ is a bounded subset of $V$ by Lemma
\ref{lemma:bounded}.  This means we can find $\eps>\delta$ such that
the ball $B_\eps$ contains $\Phi_V\inv(D)\cap S=\Phi_X\inv(D)\cap S$.
Then $\Phi_X(X\cap B_\eps)\supset\Phi_X\bigl(\Phi_X\inv(D)\cap
S\bigr)$.  By \ref{smallcone} above, $\Phi_X$ maps $S$ surjectively
onto $\Phi_X(X)$, and therefore $\Phi_X\bigl(\Phi_X\inv(D)\cap
S\bigr)=D\cap\Phi_X(X)$.  Consequently, $\Delta(X\cap B_\eps)$
contains $\Delta(X)\cap D$ and is therefore a neighbourhood of the
vertex in $\Delta(X)$.

\emph{Part\/} 2.  Suppose $\Delta(X\cap B_\eps)$ is a neighbourhood of
the vertex in $\Delta(X)$ for a certain $\eps>\delta$.  Then $\mathrm
e^{-2t}\Delta(X\cap B_\eps)$ is a neighbourhood of the vertex for all
$t$.  Choose an arbitrary $\eps'$ with $\delta<\eps'<\eps$.  By the
continuity of the retraction \eqref{equation:deformation} and the
compactness of $X\cap B_\eps$ there exists a $t$ such that $\lie
G_t(X\cap B_\eps)$ is a subset of $X\cap B_{\eps'}$.  So by
\eqref{equation:lengthflow}, $\mathrm e^{-2t}\Delta(X\cap
B_\eps)=\Delta\bigl(\lie G_t(X\cap B_\eps)\bigr)$ is a subset of
$\Delta(X\cap B_{\eps'})$, so $\Delta(X\cap B_{\eps'})$ is a
neighbourhood of the vertex in $\Delta(X)$.

\emph{Part\/} 3.  Let $D$ be any ball about the origin in $\lie k^*$.
Take any $\eps>\delta$; then there exists a $t$ such that $\mathrm
e^{-2t}\Delta(X\cap B_\eps)=\Delta\bigl(\lie G_t(X\cap B_\eps)\bigr)$
is contained in $D$.  Again by continuity and compactness, there
exists an $\eps'$ with $\delta<\eps'<\eps$ such that $X\cap
B_{\eps'}\subset\lie G_t(X\cap B_\eps)$.  But then $\Delta(X\cap
B_{\eps'})$ is contained in $D$.
\end{proof}

This proof gives no information on the shape of the set $\Delta(X\cap
B_\eps)$ away from the vertex.  It seems not unlikely that
$\Delta(X\cap B_\eps)$ is convex.  I am tempted to speculate that it is
equal up to a dilation to the momentum polytope of the projective
closure of $X$, $\Delta(\bar X)$.

\begin{example}[homogeneous vector
bundles]\label{example:associatedstein} Let $L$ be a closed subgroup
of $K$ and let $F$ be the reductive subgroup $L\co$ of $G$.  Let $W$
be a unitary $L$-module and let $X=G\times^FW$.  There exists an
orthogonal (real) representation $V_1$ of $K$ containing a vector
$v_0$ with stabilizer $K_{v_0}=L$.  The map $k\mapsto kv_0$ therefore
induces an embedding $K/L\to V_1$.  The complexification of this map
is an embedding of $K\co/L\co=G/F$ into the unitary $K$-module
$V_1\co$.  There also exists an $L$-equivariant isometric embedding of
$W$ into a unitary $K$-module $V_2$.  Then the map $X\to V_1\co\oplus
V_2$ defined by $[g,w]\mapsto gv_0+gw$ is a $G$-equivariant closed
embedding of $X$ into $V=V_1\co\oplus V_2$.  (See Lemmas 1.16 and 1.18
of \cite{sj:ho} for a proof of these facts.)  Let $\rho(v)=\|v\|^2/2$,
where $\lVert{\cdot}\rVert$ denotes the length function with respect
to the direct sum metric on $V$.  Let us identify $X$ with its image
in $V$.  Using Theorem \ref{theorem:kempfness} one can easily show
that $\rho_X$ has a unique critical level, which is a minimum, and
that the critical set is the compact orbit $Kv_0$.  Hence, by Theorem
\ref{theorem:stein}, the sets $\Delta(U)$, where $U$ ranges over the
neighbourhoods of $Kv_0$ in $X$, form a basis of neighbourhoods of the
vertex of $\Delta(X)$.
\end{example}

Here is an example of a singular variety for which the conclusion of
Theorem \ref{theorem:stein} holds.

\begin{example}\label{example:conicalfibre}
Let $F$, $W$ and $X$ be as in the previous example, and let $Y$ be an
$F$-invariant affine cone in $W$.  Let $X'$ be the affine subvariety
$G\times^FY$ of the vector bundle $X=G\times^FW$ and embed $X$ into a
$G$-module $V$ as in the previous example.  Because $Y\subset W$ is
invariant under dilations, the subvariety $X'$ is invariant under the
gradient flow of $\rho_X$.  It follows that the restriction of the flow
$\lie G_X$ to $X'$ retracts $X'$ onto the compact orbit $Kv_0$.
Exactly the same proof as that of Theorem \ref{theorem:stein} now
shows that the sets $\Delta(U')$, where $U'$ ranges over the
neighbourhoods of $Kv_0$ in $X'$, form a basis of neighbourhoods of
the vertex of the cone $\Delta(X')$.
\end{example}

\section{Hamiltonian actions and convexity}\label{section:convex}

In this section I explain how the previous, mainly algebro-geometric,
results can be generalized to arbitrary Hamiltonian actions.  The
basic idea is that every symplectic manifold with a Hamiltonian
$K$-action can locally near every orbit in the zero fibre of the
momentum map be identified with a germ of a complex affine
$G$-variety.  I then state the main result of the paper, Theorem
\ref{theorem:convex}.  First recall the following standard definition.

\begin{definition}
Let $M$ be a Hamiltonian $K$-manifold with momentum map $\Phi$.  For
every $\mu\in\lie k^*$ the (\emph{Meyer-Marsden-Weinstein\upn) reduced
space\/} or \emph{symplectic quotient\/} at level $\mu$ is the space
$\Phi\inv(K\mu)/K$.  It is denoted by $M_{\mu,K}$, or by $M_\mu$, if
the group $K$ is clear from the context.
\end{definition}

By the results of \cite{sj:st}, the symplectic quotient is a
stratified space carrying natural symplectic forms on the strata
satisfying certain compatibility conditions.  For most $\mu$ in $\lie
k^*$, $M_\mu$ is actually a symplectic V-manifold.

One application of symplectic reduction is the construction of ``local
models'', which I now briefly explain.  See \cite{ma:co} or
\cite{gu:sy} for details.  Let $\mu$ be any vector in $\lie k^*$ and
let $L$ be any closed subgroup of $K_\mu$.  Let $W$ be a symplectic
representation of $L$ and let $\Phi_W$ be the standard quadratic
$L$-momentum map on $W$.  Let $\lie z_\mu$ be the centre of $\lie
k_\mu$.  Then the zero-weight space in $\lie k$ under the adjoint
action of $\lie z_\mu$ is exactly $\lie k_\mu$, so $\lie k_\mu$ has a
natural $\lie k_\mu$-invariant complement in $\lie k$.  This means
that the principal fibre bundle $K_\mu\to K\to K\mu$ comes equipped
with a natural connection.  Furthermore, the Levi decomposition $\lie
k_\mu=\lie z_\mu\oplus[\lie k_\mu,\lie k_\mu]$ shows that $\lie z_\mu$
is a direct summand of $\lie k_\mu$.  We can therefore view $\lie
z_\mu^*$ as a subspace of $\lie k_\mu^*$ and $\lie k_\mu^*$ as a
subspace of $\lie k^*$.  Under these natural identifications, $\mu$ is
an element of $\lie z_\mu^*\subset\lie k_\mu^*$.

The manifold $K\times\lie k_\mu^*$ carries a natural closed two-form
(the minimal-coupling form defined by the symplectic form on $K\mu$
and the connection on $K\to K\mu$), which is nondegenerate in a
$K$-invariant neighbourhood of $K\times\{0\}$.  Now consider the
manifold $X=K\times\lie k_\mu^*\times T^*K_\mu\times W$ and identify
$T^*K_\mu$ with $K_\mu\times\lie k_\mu^*$ by means of
left-translations.  The action of $K_\mu\times L$ on $X$ defined by
$$
(k,l)\cdot(g,\kappa,h,\nu,w)=(gk\inv,k\kappa,khl\inv,l\nu,lw)
$$
is Hamiltonian with momentum map $\Psi\colon X\to\lie
k_\mu^*\times\lie l^*$ given by
$$
\Psi(g,\kappa,h,\nu,w)=\bigl(-\kappa+h\nu,-\nu|_{\lie
l}+\Phi_W(w)\bigr).
$$

\begin{definition}\label{definition:model}
$X(\mu,L,W)=\Psi\inv(-\mu,0)/(K_\mu\times L)$ is the symplectic
quotient of $X$ by the $K_\mu\times L$-action at the value
$(-\mu,0)\in\lie k_\mu^*\times\lie l^*$.
\end{definition}

It turns out that $X(\mu,L,W)$ is smooth.  Consider the $K$-action on
$X$ defined by left-multiplication on the first factor.  It is
Hamiltonian as well and, moreover, it commutes with the $K_\mu\times
L$-action.  It descends therefore to a Hamiltonian $K$-action on
$X(\mu,L,W)$.  The easiest way to write the $K$-momentum map on
$X(\mu,L,W)$ is as follows.  Put $\lie m=\lie k_\mu/\lie l$.  Choose
an $L$-invariant complement of $\lie l$ in $\lie k_\mu$ and identify
it with $\lie m$.  Then $\lie m^*$ is a subspace of $\lie k_\mu^*$.
The map $\varphi$ from the product $K\times\lie m^*\times W$ to
$\Psi\inv(-\mu,0)$ defined by
$$
\varphi(g,\nu,w)=\bigl(g,\nu+\Phi_W(w)+\mu,1,\nu+\Phi_W(w),w\bigr)
$$ 
descends to a $K$-equivariant diffeomorphism
$$
\bar\varphi\colon K\times^L(\lie m^*\times W)\longrightarrow
X(\mu,L,W).
$$
Via this diffeomorphism, the associated bundle $K\times^L(\lie
m^*\times W)$ acquires a closed two-form that is symplectic in a
neighbourhood of the zero section.  Note that the definition of
$\bar\varphi$ depends only on the choice of a complement of $\lie l$
in $\lie k_\mu$ (which is equivalent to the choice of a connection on
the principal fibre bundle $L\to K\to K/L$).  We shall henceforth
identify $X(\mu,L,W)$ with $K\times^L(\lie m^*\times W)$ through
$\bar\varphi$.  The $K$-momentum map on $K\times^L(\lie m^*\times W)$
is given by
\begin{equation}\label{equation:map}
\Phi([g,\nu,w])=g\bigl(\nu+\Phi_W(w)+\mu\bigr).
\end{equation}

It is useful to consider the restriction of $\Phi$ to a $K$-invariant
neighbourhood of the point $[1,0,0]$.  If we perform the reduction of
$X$ by the $K_\mu\times L$-action in stages, first with respect to
$L$ and then with respect to $K_\mu$, we see that $X(\mu,L,W)$
is an ``iterated'' associated bundle: it is a bundle
\begin{equation}\label{equation:bundle}
X(\mu,L,W)\cong K\times^{K_\mu}\bigl(K_\mu\times^L(\lie
m^*\times W)\bigr)=K\times^{K_\mu}Y,
\end{equation}
over the coadjoint orbit $K\mu\cong K/K_\mu$ with fibre the
Hamiltonian $K_\mu$-space
\begin{equation}\label{equation:fibre}
Y=K_\mu\times^L(\lie m^*\times W).
\end{equation}
The $K_\mu$-momentum map $ Y\to\lie k_\mu^*$ is the
restriction of $\Phi$ to $Y$.  We can now write $\Phi$ as the
composition of two maps,
\begin{equation}\label{equation:composition}
K\times^{K_\mu}Y \longrightarrow K\times^{K_\mu}\lie
k_\mu^*\overset{\iota}{\longrightarrow}\lie k^*,
\end{equation}
the first of which is the unique bundle map extending the
$K_\mu$-momentum map on the fibre $Y$, and the second of which is
defined by $\iota([g,\nu])=g\nu$.  

Recall that $\lie k_\mu^*$ is a \emph{slice\/} at $\mu$ for the
coadjoint action on $\lie k^*$: the restriction of $\iota$ to a
sufficiently small $K$-invariant open neighbourhood of $[1,\mu]$ is a
$K$-equivariant embedding onto an open neighbourhood of $\mu$.  This
implies that the restriction of $\Phi$ to a sufficiently small
$K$-invariant neighbourhood $U$ of $[1,0,0]$ is a bundle map of
associated bundles over $K\mu$.  Consequently, $U\cap\Phi\inv(\lie
k_\mu^*)= U\cap Y$ and the image $\Phi(U)$ is a bundle over $K\mu$
with fibre $\Phi(U\cap Y)$.

If $\mu$ happens to be in $\tplus$, then $\lie t^*$ is a subset of
$\lie k_\mu^*$, and $\tplus$ is contained in $\lie t^*_{+,\mu}$, the
positive Weyl chamber of $\lie k_\mu^*$.  In fact, $\lie
t^*_{+,\mu}\cap D=\tplus\cap D$ for a sufficiently small neighbourhood
$D$ of $\mu$ in $\lie k^*$.  It follows from this that if $U$ is small
enough
\begin{equation}\label{equation:cross1}
\Delta(U)=\Phi(U)\cap\tplus=\Phi(U\cap Y)\cap\tplus=\Phi(U\cap
Y)\cap\lie t^*_{+,\mu}=\Delta(U\cap Y),
\end{equation}
where $\Delta(U\cap Y)$ stands for the momentum set of $U\cap Y$
considered as a Hamiltonian $K_\mu$-space.

Now let $M$ be an arbitrary Hamiltonian $K$-manifold with momentum map
$\Phi$.  Marle \cite{ma:co} and Guillemin and Sternberg \cite{gu:sy}
have shown that locally at any orbit $M$ is isomorphic to some
$X(\mu,L,W)$.

\begin{theorem}[symplectic slices]\label{theorem:model}
Let $m$ be any point in $M$.  Let $L=K_m$ be the stabilizer of $m$\upn,
let $\mu=\Phi(m)$\upn, and let
$$
W=T_m(Km)^\omega\big/\bigl(T_m(Km)\cap T_m(Km)^\omega\bigr).
$$
Then there exist a $K$-invariant neighbourhood $U_1$ of $m$ in
$M$\upn, a $K$-invariant neighbourhood $U_2$ of $[1,0,0]$ in
$X(\mu,L,W)$\upn, and a map $f\colon U_1\to U_2$ with the
following properties\upn: $f$ is a $K$-equivariant
symplectomorphism\upn; $f$ intertwines the momentum maps on $U_1$ and
$U_2$\upn; and $f(m)=[1,0,0]$.
\qed
\end{theorem}

The symplectic vector space $W$ is called the \emph{symplectic
slice\/} at $m$.  Note that in the situation of the theorem $\lie m$
is simply the tangent space to $K_\mu m$ at $m$.  An immediate
consequence of the theorem is that the germ at $m$ of $\Phi\inv(\lie
k_\mu^*)$, called a \emph{local cross-section\/} of $M$, is a smooth
$K_\mu$-invariant symplectic submanifold, because in the local model
$X(\mu,L,W)$ it is equal to the germ of $Y=K_\mu\times^L(\lie
m^*\times W)$ at $[1,0,0]$.  Furthermore, if $m$ has the property that
$\mu=\Phi(m)\in\tplus$, then for any small $K$-invariant neighbourhood
$U$ of $m$ we have an equality of momentum sets
\begin{equation}\label{equation:cross2}
\Delta(U)=\Delta\bigl(U\cap\Phi\inv(\lie k_\mu^*)\bigr),
\end{equation}
because by \eqref{equation:cross1} the same is true in the local
model.  An analogue of \eqref{equation:cross2} for projective
varieties was proved by Brion (Proposition 4.1 in \cite{br:im}).

\begin{definition}\label{definition:localcone}
Let $m$ be any point in $M$.  The \emph{local momentum cone\/} at $m$
is the set $\Delta_m=\mu+\Delta(Y_m)$.  Here $Y_m$ is the complex
affine $(K_\mu)\co$-variety $(K_\mu)\co\times^{(K_m)\co}W$, with
$\mu=\Phi(m)$ and with $W$ being the symplectic slice at $m$ furnished
with a compatible $K_m$-invariant complex structure.
\end{definition}

So $\Delta_m$ is a convex polyhedral cone with vertex $\mu$ and is
contained in $\lie t^*_{+,\mu}$, the positive Weyl chamber of
$K_\mu$.  Up to a translation by $\mu$ it is a rational cone.

{}From the symplectic slice theorem we deduce the following ``local''
convexity theorem.

\begin{theorem}\label{theorem:localcone}
\begin{enumerate}
\item\label{localcone} For every $m$ in $M$ such that
$\mu=\Phi(m)$ is contained in $\tplus$ and for every sufficiently
small $K$-invariant neighbourhood $U$ of $m$ the set $\Delta(U)$ is a
neighbourhood of the vertex of the local momentum cone
$\Delta_m\subset\lie t^*_{+,\mu}$.   In particular\upn, the cone with
vertex $\mu$ spanned by $\Delta(U)$ is equal to $\Delta_m$\upn;
\item\label{samecone} for every $\mu\in\lie k^*$ and for every
connected component $C$ of the fibre $\Phi\inv(\mu)$\upn, the local
momentum cone $\Delta_m$ is independent of the point $m\in C$.
\end{enumerate}
\end{theorem}

\begin{proof}
\ref{localcone}. First we treat the case $\mu=0$.  At points in the
zero fibre of $\Phi$ there is an alternative local model for $M$ as
follows.  Put $F=L\co=(K_m)\co$ and choose an $F$-invariant compatible
complex structure on the symplectic slice $W$ at $m$.  Consider the
affine $G$-variety $Z=G\times^FW$.  Embed it into a unitary $K$-module
$V$ as in Example \ref{example:associatedstein}, and regard it as a
Hamiltonian $K$-space with the symplectic structure and momentum map
inherited from $V$.  The stabilizer of the point $[1,0]\in Z$ under
the $K$-action is equal to $L$.  The $G$-orbit through $[1,0]$ is
closed, so by Theorem \ref{theorem:kempfness}, $[1,0]$ is contained in
$\Phi_Z\inv(0)$.  The symplectic slice at $[1,0]$ is simply the
(Hermitian) orthogonal complement to the $G$-orbit through $[1,0]$ and
is therefore equal to $W$.  Theorem \ref{theorem:model} now allows us
to conclude that the $K$-invariant germ of $M$ at $m$ is isomorphic as
a Hamiltonian $K$-manifold to the $K$-invariant germ of $Z$ at
$[1,0]$.

Putting this information together with Theorems \ref{theorem:affine}
and \ref{theorem:stein} and Example \ref{example:associatedstein}, we
see that for every small $K$-invariant neighbourhood $U$ of $m$ the
set $\Delta(U)=\Phi(U)\cap\tplus$ is a neighbourhood of the vertex of
the cone $\Delta_m=\Delta(G\times^FW)$.

If $\mu\neq0$, we may without loss of generality assume that
$\mu\in\tplus$, because every $K$-orbit in $\lie k^*$ intersects
$\tplus$.  This case can then be reduced to the case $\mu=0$ by using
\eqref{equation:cross2} and shifting the $K_\mu$-momentum map on the
space $Y$ in \eqref{equation:fibre} by the vector $-\mu$.  This
shifted map is still an equivariant $K_\mu$-momentum map, because
$\mu$ is in $\lie z_\mu^*$.

\ref{samecone}. Again, it suffices to prove this for $\mu=0$.  If $M$
is an affine $G$-variety, then by Theorem \ref{theorem:kempfness} the
points in $\Phi\inv(0)$ are exactly those whose $G$-orbits are closed.
Corollary \ref{corollary:samecone} then says that the local momentum
cones $\Delta_m$ are constant along the fibre $\Phi\inv(0)$.  For
general $M$, this argument combined with the symplectic slice theorem
shows that $\Delta_m$ is locally constant along the fibre
$\Phi\inv(0)$.
\end{proof}

\begin{remark}
It follows from the Cartan decomposition of $G$ that there exists a
global $K$-equivariant diffeomorphism between the two local models
$K\times^L(\lie m^*\times W)$ and $G\times^FW$, but I don't know if
there is a global symplectomorphism.
\end{remark}

The following generalization of Kirwan's convexity theorem
\cite{ki:con} is the main result of this paper.  The first part
describes $\Delta(M)$ as a locally finite intersection of polyhedral
cones, each of which is determined by local data on $M$.  (The fact
that $\Delta(M)$ is a convex set when $\Phi$ is proper was also proved
in \cite{hi:sy}.)  The second part is a necessary condition for a
point to be a vertex of $\Delta(M)$.  It generalizes the well-known
fact that for torus actions vertices arise as images of fixed points.
The third part states that the points where this necessary condition
is fulfilled form a discrete subset of $\tplus$.

\begin{theorem}\label{theorem:convex}
Assume that the momentum map $\Phi\colon M\to\lie k^*$ is proper.
\begin{enumerate}
\item\label{convex}
$\Delta(M)$ is the intersection of local momentum cones\upn:
\begin{equation}\label{equation:intersection}
\Delta(M)=\bigcap_{m\in\Phi\inv(\tplus)}\Delta_m.
\end{equation}
This intersection is locally finite and therefore $\Delta(M)$ is a
closed convex polyhedral subset of $\tplus$\upn;
\item\label{extreme}
if $\mu$ is a vertex of $\Delta(M)$ and $m$ is any point in the fibre
$\Phi\inv(\mu)$\upn, then $\lie k_\mu =[\lie k_\mu,\lie k_\mu]+\lie
k_m$\upn, or\upn, equivalently\upn, $K_\mu=[K_\mu,K_\mu]K_m$.  In
particular\upn, if $\mu$ is a vertex of $\Delta(M)$ lying in the
interior of $\tplus$\upn, then $T$ fixes $m$\upn;
\item\label{compact}
let $E$ be the subset of $M$ consisting of all points $m$ such that
$\mu\in\tplus$ and $\lie k_\mu=[\lie k_\mu,\lie k_\mu]+\lie k_m$\upn,
where $\mu=\Phi(m)$.  The image $\Phi(E)$ is a discrete subset of
$\tplus$.  If $M$ is compact\upn, then $\Delta(M)$ is the convex hull
of $\Phi(E)$.
\end{enumerate}
\end{theorem}

\begin{proof}
\ref{convex}.  The assumption that $\Phi$ is proper implies that its
image $\Phi(M)$ is closed and, by the argument outlined in Section
\ref{subsection:flow}, that its fibres are connected.  Consequently,
by Theorem \ref{theorem:localcone}.\ref{samecone}, for every
$\mu\in\Delta(M)$ the cone $\Delta_m$ is the same for all points
$m\in\Phi\inv(\mu)$.  It is now easy to deduce from Theorem
\ref{theorem:localcone}.\ref{localcone} plus the fact that $\Phi$ is
proper that for every $\mu\in\Delta(M)$ there exists an open subset
$D$ of $\lie t^*$ containing $\mu$ such that
\begin{equation}\label{equation:locallyconvex}
\Delta(M)\cap D=\Delta_m\cap D,
\end{equation}
where $m$ is any point in the fibre $\Phi\inv(\mu)$.  This means that
$\Delta(M)$ is locally convex.  But every closed locally convex set is
convex, so $\Delta(M)$ is convex.  Since every closed convex set is
the intersection of all closed cones containing it, the equality
\eqref{equation:locallyconvex} also implies
\eqref{equation:intersection}.  Furthermore, applying Theorem
\ref{theorem:localcone} to the Hamiltonian $K$-manifold $\Phi\inv(D)$,
we find that for every $\mu'\in D$ and every $m'\in\Phi\inv(\mu')$ the
local momentum cone $\Delta_{m'}$ is equal to the cone with vertex
$\mu'$ spanned by $\Delta_m$.  Since $\Delta_m$ is a polyhedral cone,
it follows from this that as $\mu'$ ranges over $D$, only finitely
many different cones $\Delta_{m'}$ can occur.  In other words, the
collection of cones appearing in the intersection
\eqref{equation:intersection} is locally finite on $\tplus$.  This
means that $\Delta(M)$ is a polyhedron.

\ref{extreme}.  This follows from \eqref{equation:locallyconvex} and
Theorem \ref{theorem:propercone} (applied to the group $G=(K_\mu)\co$
and the variety $X=(K_\mu)\co\times^{(K_m)\co}W$).

\ref{compact}.  First I prove that $\Phi(E)$ is discrete.  Since
$\Phi$ is proper, it suffices to show that every point $m$ in $M$
possesses a neighbourhood $U$ such that $\Phi(E\cap U)$ is discrete.
By the symplectic slice theorem, we may therefore assume $M$ is an
affine variety.  But it follows from Lemma \ref{lemma:extreme} below
that for an affine $G$-variety the set $\Phi(E)$ consists of the
origin in $\tplus$ only.

Finally, the second statement in \ref{compact} follows immediately
from \ref{extreme} and the fact that a compact convex set is the
convex hull of its extreme points.
\end{proof}

\begin{lemma}\label{lemma:extreme}
Let $V$ be a unitary $K$-module and suppose that $v\in V$ satisfies
the condition $K_\mu=[K_\mu,K_\mu]K_v$\upn, where $\mu=\Phi_V(v)$.
Then $\mu=0$.
\end{lemma}

\begin{proof}
We want to prove that $\mu(\xi)=\Phi_V^\xi(v)=(\sq/2)\langle\xi
v,v\rangle=0$ for all $\xi\in\lie k$.  It suffices to show this for
all $\xi\in\lie k_\mu$.  The condition on $v$ says that
$\xi=[\chi,\zeta]+\eta$ for some $\chi,\zeta\in[\lie k_\mu,\lie
k_\mu]$ and $\eta\in\lie k_v$.  Now $\eta v=0$, so $\mu(\eta)=0$; and
$\mu(\chi,\zeta)=\mu(\ad\chi(\zeta))=(\ad^*\chi)\mu(\zeta)=0$.  We
conclude that $\mu(\xi)=0$.
\end{proof}

Theorem \ref{theorem:convex} is clearly not optimal.  As Theorem
\ref{theorem:affine} shows, it is not always necessary to assume that
$\Phi$ is proper, nor even that $M$ is nonsingular.  See Section
\ref{section:examples} for further examples of convexity for
noncompact or singular spaces.

On the other hand, the necessary condition \ref{extreme} for a point
to be mapped to a vertex of $\Delta(M)$ is optimal in the following
sense.

\begin{proposition}\label{proposition:vertex}
For every $\mu\in\tplus$ and for every closed subgroup $L$ of $K_\mu$
such that $K_\mu=[K_\mu,K_\mu]L$\upn, there exists a Hamiltonian
$K$-manifold $(M,\omega,K,\Phi)$ with a point $m\in M$ satisfying the
following properties\upn: $K_m=L$\upn, $\Phi(m)=\mu$\upn, and
$\lambda$ is a vertex of $\Delta(M)$.
\end{proposition}

The proof will be given in Section \ref{subsection:reduction}.

Theorem \ref{theorem:convex}.\ref{convex} implies that for every $\mu$
in $\Delta(M)$ the cone with vertex $\mu$ spanned by $\Delta(M)$ is
equal to the local momentum cone $\Delta_m$, where $m$ is any point in
$\Phi\inv(\mu)$.  In other words, the local shape of $\Delta(M)$ near
$\mu$ is determined by the representation of the isotropy subgroup
$K_m$ on the symplectic slice $W$ at $m$.  Calculating $\Delta_m$
boils down to finding generators for the monoid of highest weights of
the homogeneous vector bundle $(K_\mu)\co\times^{(K_m)\co}W$.  If
$\mu$ is contained in the boundary of the positive Weyl chamber, this
is usually an arduous task.  However, the situation is more manageable
if the fibre $\Phi\inv(\mu)$ contains a point $m$ that is fixed under
$K_\mu$.  This means that $K_m=K_\mu$, or, equivalently, that the
restriction of $\Phi$ to the orbit $Km$ is a symplectic isomorphism
onto the coadjoint orbit $K\mu$.  If this is the case, the vector
space $\lie m$ in the local model \eqref{equation:bundle} is $0$, so
the tangent space at $m$ to the symplectic cross-section
$\Phi\inv(\lie k_\mu^*)$ is equal to the symplectic slice $W$ at $m$,
and $W$ is simply the symplectic orthogonal complement of $T_m(Km)$
inside $T_mM$.  In other words,
\begin{equation}\label{equation:symplecticorbit}
T_m\bigl(\Phi\inv(\lie k_\mu^*)\bigr)=W=T_mM/T_m(Km)\cong
\bigl(T_m(Km)\bigr)^\omega.
\end{equation}
It follows that $\Delta_m=\Delta(W)$, where $W$ is regarded as a
$K_\mu$-module.  In some examples this enables one to determine the
entire momentum set $\Delta(M)$; see for instance Section
\ref{subsection:projective}.

\section{Examples}\label{section:examples}

\subsection{Actions on projective space}
\label{subsection:projective}

Let $V$ be a finite-dimensional unitary $K$-module.  If $V$ is
irreducible and has highest weight $\lambda$, it follows from
\eqref{equation:weyl} that the momentum polytope of the projective
space $\P V$ for the $T$-action is the convex hull of the Weyl group
orbit through $\lambda^*$.  Similarly, if $\lambda_1$,
$\lambda_2,\dots$, $\lambda_k$ are the highest weights of the
irreducible submodules of $V$, then the $T$-momentum map image of $\P
V$ is the convex hull of the union of the $\lie W$-orbits through
$\lambda_1^*$, $\lambda_2^*,\dots$, $\lambda_k^*$.  This implies
$\Delta(\P V)$ is a subset of
$$
\tplus\cap\hull(\lie W\lambda_1^*\cup\lie
W\lambda_2^*\cup\cdots\cup\lie W\lambda_k^*).
$$
Arnal and Ludwig \cite{ar:co} and Wildberger \cite{wi:mo} determined
this subset for ``most'' $V$ that are irreducible.  I shall calculate
$\Delta(\P V)$ in some (but not all) of the remaining cases.

For simplicity I assume $K$ to be semisimple, although the results
can easily be generalized to arbitrary compact groups.  Let $\Psi$ be
the root system of $(\lie k,\lie t)$ and let $\lie
g=\bigoplus_{\alpha\in\Psi}\C E_\alpha$ be the root space
decomposition of $\lie g=\lie k\co$.  The following result can be used
to find a ``lower bound'' for the polytope $\Delta(\P V)$.

\begin{lemma}\label{lemma:unitary}
\begin{enumerate}
\item\label{vector} Let $v$ be any vector in $V$ and let $[v]$ be the
ray through $v$.  Then $\Phi_{\P V}([v])\in\lie t^*$ if and only
if $\langle E_\alpha v,v\rangle=0$ for all roots $\alpha$.
\item\label{subspace} Let $v_1,\dots$\upn, $v_l$ be weight vectors in
$V$ with weights $\nu_1,\dots$\upn, $\nu_l$.  Assume that $\langle
E_\alpha v_i,v_j\rangle=0$ for all roots $\alpha$ and for all $i$ and
$j$ with $1\leq i<j\leq l$.  \upn(This is for instance the case if for
all $i$ and $j$ the difference $\nu_i-\nu_j$ is not a root.\upn) Then
the subspace spanned by $v_1,\dots$\upn, $v_l$ is contained in
$\Phi_{V}\inv(\lie t^*)$\upn, and the intersection
$\tplus\cap\hull\{\nu_1^*,\dots,\nu_l^*\}$ is contained in $\Delta(\P
V)$.
\end{enumerate}
\end{lemma}

\begin{proof}
\ref{vector}.  Let $\xi_\alpha=E_\alpha-E_{-\alpha}$ and
$\eta_\alpha=\sq(E_\alpha+E_{-\alpha})$.  By
\eqref{equation:projective1}, $\Phi_{\P V}([v])\in\lie t^*$ if and
only if $\langle\xi_\alpha v,v\rangle =\langle\eta_\alpha
v,v\rangle=0$ for all roots $\alpha$.  Since $E_{-\alpha}$ (viewed as
an operator on $V$) is the adjoint of $E_\alpha$, this is equivalent
to $\langle E_\alpha v,v\rangle=0$ for all $\alpha$.

\ref{subspace}.  Let $W$ denote the linear span of the $v_i$.  Let
$\alpha$ be any root.  The assumption implies $\langle E_\alpha
v,v\rangle=0$ whenever $v$ is in $W$.  It now follows from
\ref{vector} that $\Phi_V(W)$ is contained in $\lie t^*$.

{}From this fact and from \eqref{equation:projective1} we infer that
$$
\Phi_{\P V}(g[c_1v_1+\cdots+c_lv_l])
=\frac{|c_1|^2\nu_1^*+\cdots+|c_l|^2\nu_l^*}{|c_1|^2+\cdots+|c_l|^2},
$$
where $g$ is any element of the normalizer of $T$ representing
$w_0\in\eu W$.  Hence $\Delta(\P V)=\Phi_{\P V}(\P V)\cap\tplus$
contains the set $\tplus\cap\hull\{\nu_1^*,\dots,\nu_l^*\}$.
\end{proof}

Put $\check\mu=2\mu/(\mu,\mu)$ for any $\mu$ in $(\lie t\co)^*$.  Let
$\alpha_1,\dots$, $\alpha_r$ be the simple roots of $\lie k$ and
$\pi_1,\dots$, $\pi_r$ the corresponding fundamental weights, that is,
the basis of $\lie t^*$ dual to $\check\alpha_1,\dots$,
$\check\alpha_r$.  Then
$\lambda=\sum_{i=1}^r(\lambda,\check\alpha_i)\pi_i$, where the
coefficients $(\lambda,\check\alpha_i)$ are nonnegative integers.  The
next result gives an upper bound for the polytope $\Delta(\P V)$
for irreducible $V$.

\begin{proposition}\label{proposition:unirrep}
Assume $V$ is irreducible and has highest weight $\lambda$.  Let
$\Pi_\lambda$ be the set of weights of $V$ that are not of the form
$\lambda-\alpha$ for any positive root $\alpha$ such that
$(\lambda,\check\alpha)=1$.  Then $\Delta(\P V)$ is contained in
$\tplus\cap\hull\Pi_\lambda^*$.
\end{proposition}

\begin{proof}
Let $v_\lambda$ be a highest-weight vector in $V$, so that $\Phi_{\P
V}([v_\lambda])=-\lambda$.  It suffices to show that the local
momentum cone $\Delta_{[v_\lambda]}$ of $\P V$ at the point
$[v_\lambda]$ (see Definition \ref{definition:localcone}) is contained
in the cone with vertex $-\lambda$ spanned by the set $-\Pi_\lambda$.
To this end let us compute the tangent space to the symplectic
cross-section at $[v_\lambda]$ and the weights of the $T$-action on
it.  The vector $v_\lambda$ is an eigenvector for $K_\lambda$, so
$K_{[v_\lambda]}=K_\lambda$.  In view of
\eqref{equation:symplecticorbit} this implies that we have an
isomorphism of $K_\lambda$-modules
\begin{equation}\label{equation:slice}
T_{[v_\lambda]}\bigl(\Phi_{\P V}\inv(\lie k_\lambda^*)\bigr)=W \cong
T_{[v_\lambda]}(\P V)/T_{[v_\lambda]}(K[v_\lambda]),
\end{equation}
$W$ being the symplectic slice at $[v_\lambda]$.  Define $\Pi$ to be
the subset of the weight lattice $\Lambda^*$ consisting of $0$ and of
all weights occurring in $W$ under the action of the maximal torus $T$
of $K_\lambda$.  From \eqref{equation:slice} we get:
\begin{equation}\label{equation:localcone1}
\Delta_{[v_\lambda]}\subset-(\lambda+\cone\Pi).
\end{equation}
(If $\lambda$ is strictly dominant, so that $K_\lambda=T$, then this
inclusion is an equality.)  I assert that
\begin{equation}\label{equation:localcone2}
\cone\Pi=\cone(-\lambda+\Pi_\lambda).
\end{equation}
Because of \eqref{equation:localcone1}, this will finish the proof.
Let $V=\bigoplus_{\nu\in\Lambda^*}V_\nu$ be the weight space
decomposition of $V$ and let $\C_{-\lambda}$ be the one-dimensional
representation of $K_\lambda$ defined by the character
$-\lambda\in\lie z_\lambda^*$.  Then the quotient map $V-\{0\}\to\P V$
induces an isomorphism of $K_\lambda$-modules
\begin{equation}\label{equation:tangent1}
T_{[v_\lambda]}(\P V) \cong
\bigoplus_{\nu\in\Lambda^*-\{\lambda\}}V_\nu\otimes\C_{-\lambda}.
\end{equation}
Since the maximal unipotent subgroup $N$ fixes $v_\lambda$, the
complex stabilizer $G_{[v_\lambda]}$ is the parabolic subgroup
$P_\lambda=(K_\lambda)\co N$.  This implies that the real orbit
$K[v_\lambda]$ is equal to the complex orbit $G[v_\lambda]$ and
therefore we have natural isomorphisms of complex
$K_\lambda$-representations
\begin{equation}\label{equation:tangent2}
T_{[v_\lambda]}(K[v_\lambda])\cong\lie p_\lambda^\circ
\cong\bigoplus_{\substack{\alpha\in\Psi^-\\ (\lambda,\alpha)\neq0}} \C
E_{\alpha},
\end{equation}
where $\lie p_\lambda^\circ$ denotes the annihilator of $\lie
p_\lambda$ in $\lie g^*$.  (The second isomorphism is induced by the
Killing form on $\lie g$.)  Let $\Lambda^*_{\text r}\subset\Lambda^*$
denote the root lattice of the pair $(\lie k,\lie t)$.  There are two
inclusions,
$$
(\Lambda^*_{\text r}\backslash\Psi^-)\cap(-\lambda+\hull\lie
W\lambda)\subset\Pi\subset\Lambda^*_{\text r}\cap(-\lambda+\hull\lie
W\lambda),
$$
the second of which follows from \eqref{equation:tangent1} and the
first of which follows from \eqref{equation:slice} and
\eqref{equation:tangent2}.  Moreover, by the definition of
$\Pi_\lambda$, the same inclusions hold with $\Pi$ replaced by
$-\lambda+\Pi_\lambda$.  Therefore, to establish
\eqref{equation:localcone2}, it suffices to prove the following
statement for every positive root $\alpha$: if $-\alpha\in\Pi$, then
$-\alpha\in-\lambda+\Pi_\lambda$; and if
$-\alpha\in-\lambda+\Pi_\lambda$, then some multiple of $-\alpha$ is
in $\Pi$.  There are three possibilities:

1.  $(\lambda,\check\alpha)=0$.  Then $\lambda-\alpha$ is not a weight
of $V$, so $-\alpha$ is contained in neither $\Pi$ nor
$-\lambda+\Pi_\lambda$.

2.  $(\lambda,\check\alpha)=1$.  Then
$\lambda-\alpha\not\in\Pi_\lambda$ by definition.  Also,
$\lambda-\alpha$ is a weight of $V$ with multiplicity one, so by
\eqref{equation:tangent2} $-\alpha$ is not a weight of $W$, so
$-\alpha\not\in\Pi$.

3.  $(\lambda,\check\alpha)=l>1$.  Then the weights $\lambda$,
$\lambda-\alpha,\dots$, $\lambda-l\alpha$ occur in $V$.  Since
$(\lambda,\check\alpha)\neq1$, all these weights are also in
$\Pi_\lambda$.  The weight $-\alpha$ may or may not occur in $W$
(depending on whether the multiplicity of $\lambda-\alpha$ in $V$ is
greater than one), but at any rate $-l\alpha$ is a weight of $W$.  We
conclude that $-\alpha\in-\lambda+\Pi_\lambda$ and $-l\alpha\in\Pi$.

In sum, we have shown that $\Pi\subset-\lambda+\Pi_\lambda$ and
$-\lambda+\Pi_\lambda\subset\cone\Pi$.  This proves
\eqref{equation:localcone2}.
\end{proof}

The following result follows easily from Lemma \ref{lemma:unitary} and
Proposition \ref{proposition:unirrep}.

\begin{proposition}[\cite{ar:co},\cite{wi:mo}]
\label{proposition:bigweight}
Assume $V$ is irreducible and has highest weight $\lambda$.
\begin{enumerate}
\item Suppose that $(\lambda,\check\alpha_i)\neq1$ for $i=1$\upn,
$2,\dots$\upn, $r$.  Then $\Delta(\P V)=\tplus\cap\hull\lie
W\lambda^*$\upn;
\item if $(\lambda,\check\alpha)=1$ for some positive root
$\alpha$\upn, then $\Delta(\P V)$ is not equal to $\tplus\cap\hull\lie
W\lambda^*$.
\item $\Phi_{\P V}(\P V)\cap\lie t^*$ is convex if and only if
$(\lambda,\check\alpha_i)\neq1$ for $i=1$\upn, $2,\dots$\upn, $r$.
\qed
\end{enumerate}
\end{proposition}



If $(\lambda,\check\alpha)=1$ for some positive root $\alpha$,
Proposition \ref{proposition:bigweight} does not give an explicit
description of the polytope $\Delta(\P V)$.  Using Lemma
\ref{lemma:unitary} one can easily check that the upper bound given by
Proposition \ref{proposition:unirrep} is sharp e.g.\ in the following
cases.

\begin{proposition}\label{proposition:example}
The equality $\Delta(\P V)=\tplus\cap\hull\Pi_\lambda^*$ holds
if
\begin{enumerate}
\item $K=\SU(4)$ and $\lambda=\pi_1$\upn, $\pi_2$\upn, $\pi_3$\upn, or
$\pi_1+\pi_2+\pi_3$\upn;
\item $K$ has rank two and $\lambda\in\tplus$ is arbitrary.
\qed 
\end{enumerate}
\end{proposition}

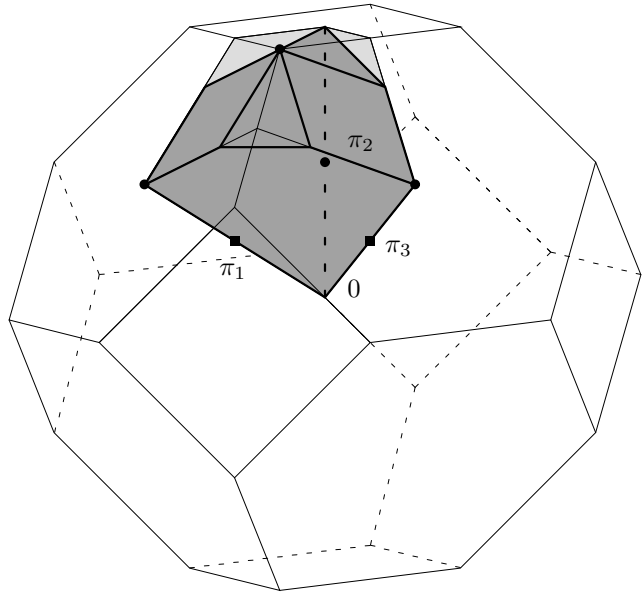
\begin{figure}
\setlength{\unitlength}{0.2mm}
$$
\begin{picture}(420,405)
\path(150,255)(60,165)(150,75)
	(240,165)(150,255)
\dashline{4.000}(270,315)(180,225)(270,135)
	(360,225)(270,315)
\dashline{4.000}(270,315)(360,225)
\path(150,255)(180,360)(120,375)
	(30,285)(0,180)(60,165)
\path(0,180)(30,105)(120,15)
	(180,0)(150,75)
\path(240,165)(360,180)(390,105)
	(300,15)(180,0)
\path(120,375)(240,390)(300,375)(180,360)
\path(300,375)(390,285)(360,180)
\path(390,285)(420,210)(390,105)
\dashline{4.000}(240,390)(270,315)
\dashline{4.000}(360,225)(420,210)
\dashline{4.000}(120,15)(240,30)(300,15)
\dashline{4.000}(270,135)(240,30)
\dashline{4.000}(180,225)(60,210)
\dashline{4.000}(30,285)(60,210)(30,105)
\texture{c0c0c0c0 0 0 0 0 0 0 0 
	c0c0c0c0 0 0 0 0 0 0 0 
	c0c0c0c0 0 0 0 0 0 0 0 
	c0c0c0c0 0 0 0 0 0 0 0 }
\shade\path(210,375)(150,367.5)(90,270)(210,195)(270,270)(240,367.5)
(210,375)
\path(210,375)(150,367.5)(90,270)(210,195)(270,270)(240,367.5)(210,375)
\path(270,270)(165,307.5)(90,270)
\texture{cccccccc 0 0 0 cccccccc 0 0 0 
	cccccccc 0 0 0 cccccccc 0 0 0 
	cccccccc 0 0 0 cccccccc 0 0 0 
	cccccccc 0 0 0 cccccccc 0 0 0 }
\put(210,195){\makebox(200,200)[lb]{\shade\path(0,0)(60,75)(40,140)%
(0,180)(-30,165)(-80,140)(-120,75)(0,0)}}
\thicklines
\put(210,195){\makebox(200,200)[lb]{\path(0,0)(60,75)(40,140)(0,180)%
(-30,165)(-80,140)(-120,75)(0,0)}}
\put(210,195){\makebox(200,200)[lb]{\path(-120,75)(-70,100)(-10,100)%
(60,75)}}
\put(210,195){\makebox(200,200)[lb]{\path(-10,100)(-30,165)(-70,100)}}
\put(210,195){\makebox(200,200)[lb]{\path(-30,165)(40,140)}}
\put(150,232.5){\makebox(0,0)[0mm]{\rule{1.2mm}{1.2mm}}} 
\put(210,285){\circle*{5}} 
\put(240,232.5){\makebox(0,0)[0mm]{\rule{1.2mm}{1.2mm}}} 
\put(180,360){\circle*{5}} 
\put(90,270){\circle*{5}} 
\put(270,270){\circle*{5}} 
\put(225,195){\makebox(0,0)[lb]{\smash{$0$}}}
\put(140,210){\makebox(0,0)[lb]{\smash{$\pi_1$}}}
\put(225,295){\makebox(0,0)[lb]{\smash{$\pi_2$}}}
\put(250,225){\makebox(0,0)[lb]{\smash{$\pi_3$}}}
\dashline{6.000}(210,195)(210,375)
\thinlines
\path(60,165)(150,255)(180,360)
\path(150,367.5)(180,360)(240,367.5)
\path(150,255)(240,165)
\path(270,270)(165,307.5)(90,270)
\end{picture}
$$
\caption{$K=\protect\operatorname{SU}(4)$ and $\lambda=\pi_1+\pi_2+\pi_3
=\frac{1}{2}\sum_{\alpha\in\Psi^+}\alpha$.   Vertices of
$\Delta(\P V)$ are $\pi_1+\pi_2+\pi_3$, $0$, $2\pi_1$, $2\pi_2$,
$2\pi_3$, $\frac{4}{3}\pi_1+\pi_2$, $\frac{4}{3}\pi_3+\pi_2$,
$\pi_1+\frac{5}{3}\pi_3$ and $\frac{5}{3}\pi_1+\pi_3$ }
\label{diagram:su4}
\end{figure}

Diagrams \ref{diagram:su3}--\ref{diagram:su4} illustrate Proposition
\ref{proposition:example}.  The figures show the convex hull of $\lie
W\lambda^*$ and the polytope $\Delta(\P V)$ (shaded).  The
intersection $\tplus\cap\hull\lie W\lambda^*$ is indicated in light
shading.  The dominant weights occurring in $V$ are denoted by black
circles.  These are exactly the images of the $T$-fixed points in $\P
V$ intersected with $\tplus$.  Notice that few of the vertices on the
walls of $\tplus$ arise as images of $T$-fixed points.  In those cases
where they are not weights of $V$, the fundamental weights of $\lie k$
are indicated by black squares.

\begin{figure}
\setlength{\unitlength}{0.1mm}
$$
\begin{picture}(600,600)(-300,-300)
\path(50,259.80762)(-250,86.60254)(-250,-86.60254)(50,-259.80762)%
(200,-173.20508)(200,173.20508)(50,259.80762)
\dashline{5}(150,259.80762)(-150,-259.80762)
\dashline{5}(-300,0)(300,0)
\dashline{5}(-150,259.80762)(150,-259.80762)
\texture{c0c0c0c0 0 0 0 0 0 0 0 
	c0c0c0c0 0 0 0 0 0 0 0 
	c0c0c0c0 0 0 0 0 0 0 0 
	c0c0c0c0 0 0 0 0 0 0 0 }
\shade\path(50,259.80762)(-100,173.20508)(0,0)(125,216.50635)%
(50,259.80762)
\path(50,259.80762)(-100,173.20508)(0,0)(125,216.50635)%
(50,259.80762)
\texture{cccccccc 0 0 0 cccccccc 0 0 0 
	cccccccc 0 0 0 cccccccc 0 0 0 
	cccccccc 0 0 0 cccccccc 0 0 0 
	cccccccc 0 0 0 cccccccc 0 0 0 }
\shade\path(50,259.80762)(-100,173.20508)(0,0)(100,173.20508)%
(50,259.80762)
\thicklines
\shade\path(50,259.80762)(-100,173.20508)(0,0)(100,173.20508)%
(50,259.80762)
\put(50,86.60254){\circle*{10}}
\put(-50,86.60254){\makebox(0,0)[0mm]{\rule{1.2mm}{1.2mm}}}
\put(50,259.80762){\circle*{10}}
\put(-100,173.20508){\circle*{10}}
\put(60,80){\makebox(0,0)[lb]{\smash{$\pi_2$}}}
\put(-100,80){\makebox(0,0)[lb]{\smash{$\pi_1$}}}
\put(65,265){\makebox(0,0)[lb]{\smash{$\lambda^*=\pi_1+2\pi_2$}}}
\thinlines
\dottedline{8}(100,173.20508)(200,0)
\end{picture}
$$
\caption{$K=\protect\operatorname{SU}(3)$ and $\lambda=2\pi_1+\pi_2$}
\label{diagram:su3}
\end{figure}

\begin{figure}
\setlength{\unitlength}{0.15mm}
$$
\begin{picture}(400,400)(-200,-200)
\path(0,173.20508)(-150,86.60254)(-150,-86.60254)(0,-173.20508)%
(150,-86.60254)(150,86.60254)(0,173.20508)
\dashline{5}(100,173.20508)(-100,-173.20508)
\dashline{5}(-200,0)(200,0)
\dashline{5}(100,-173.20508)(-100,173.20508)
\dashline{5}(0,-173.20508)(0,173.20508)
\dashline{5}(-150,86.60254)(150,-86.60254)
\dashline{5}(-150,-86.60254)(150,86.60254)
\texture{c0c0c0c0 0 0 0 0 0 0 0 
	c0c0c0c0 0 0 0 0 0 0 0 
	c0c0c0c0 0 0 0 0 0 0 0 
	c0c0c0c0 0 0 0 0 0 0 0 }
\shade\path(0,173.20508)(0,0)(75,129.90381)(0,173.20508)
\path(0,173.20508)(0,0)(75,129.90381)(0,173.20508)
\texture{cccccccc 0 0 0 cccccccc 0 0 0 
	cccccccc 0 0 0 cccccccc 0 0 0 
	cccccccc 0 0 0 cccccccc 0 0 0 
	cccccccc 0 0 0 cccccccc 0 0 0 }
\shade\path(0,173.20508)(0,0)(50,86.60254)(0,173.20508)
\thicklines
\shade\path(0,173.20508)(0,0)(50,86.60254)(0,173.20508)
\put(0,0){\circle*{7}}
\put(50,86.60254){\circle*{7}}
\put(0,173.20508){\circle*{7}}
\put(60,80){\makebox(0,0)[lb]{\smash{$\pi_1$}}}
\put(10,180){\makebox(0,0)[lb]{\smash{$\pi_2$}}}
\thinlines
\dottedline{8}(50,86.60254)(150,-86.60254)
\end{picture}
$$
\caption{$K=\protect\operatorname{G}_2$ and $\lambda=\pi_2$ (highest
weight of complexified adjoint representation, $\lie g_2^{\C}$)}
\label{diagram:g2}
\end{figure}

Finally, here is a generalization of Proposition
\ref{proposition:bigweight} to reducible representations.

\begin{proposition}
\begin{enumerate}
\item\label{bigpolytope} Let $\lambda_1$, $\lambda_2,\dots$,
$\lambda_k$ be the highest weights of the irreducible submodules of
$V$.  Suppose that $(\lambda_j,\check\alpha_i)\neq1$ for all $i$ and
$j$.  Then $\Delta(\P V)=\tplus\cap\hull(\lie W\lambda_1^*\cup\lie
W\lambda_2^*\cup\cdots\cup\lie W\lambda_k^*)$.
\item\label{multiple} Suppose $V$ is the direct sum of at least $r+1$
copies of a unitary irreducible representation with highest weight
$\lambda$\upn, where $r$ is the rank of $K$.  Then $\Delta(\P
V)=\tplus\cap\hull\lie W\lambda^*$.
\end{enumerate}
\end{proposition}

\begin{proof}
\ref{bigpolytope}.  This is easy to deduce from Lemma
\ref{lemma:unitary} and Proposition \ref{proposition:bigweight}.

\ref{multiple}.  Write $V=\bigoplus_1^kV'$, where $k>r$ and $V'$ is an
irreducible module with highest weight $\lambda$.  Clearly, $\Delta(\P
V)$ is a subset of $\tplus\cap\hull\lie W\lambda^*$.  The weight
polytope $\hull\lie W\lambda$ has exactly $r$ edges containing the
vertex $\lambda$.  Let $\nu_1$, $\nu_2,\dots$, $\nu_r$ be the opposite
endpoints of these edges.  Choose weight vectors $v_0$, $v_1,\dots$,
$v_r$, each coming from a different copy of $V'$, and having weights
$\lambda$, $\nu_1$, $\nu_2,\dots$, $\nu_r$, respectively.  Since each
copy of $V'$ is $K$-invariant and they are all mutually orthogonal,
$\langle E_\alpha v_i,v_j\rangle=0$ for all roots $\alpha$ and for all
$i$ and $j$ with $0\leq i<j\leq r$.  Lemma \ref{lemma:unitary}
therefore tells us that
$\tplus\cap\hull\{\lambda^*,\nu_1^*,\dots,\nu_r^*\}$ is contained in
$\Delta(\P V)$.  Hence, $\Delta(\P V)=\tplus\cap\hull\lie W\lambda^*$.
\end{proof}

\subsection{Cotangent bundles}\label{subsection:cotangent}

Let $Q$ be a connected $K$-manifold and let $M$ be the cotangent
bundle of $Q$.  Points in $M$ will be written as pairs $(q,p)$, where
$q\in Q$ and $p\in T^*_qQ$, and tangent vectors to $M$ as pairs
$(\delta q,\delta p)$, where $\delta q\in T_qQ$ and $\delta p\in
T_p(T^*_qQ)$.  The standard one-form $\alpha$ on $M$ is the
$K$-invariant form defined by $\alpha_{(q,p)}(\delta q,\delta
p)=p(\delta q)$.  The two-form $\omega=-d\alpha$ is symplectic, and
the lifted $K$-action on $M$ is Hamiltonian with momentum map defined
by $\Phi^\xi =\iota_{\xi_M}\alpha$, that is,
$\Phi^\xi(q,p)=p(\xi_{Q,q})$.  Clearly, $\Phi$ is homogeneous of
degree one in $p$, so $\Phi\inv(0)$ is a conical subset of $M$.  In
particular, $\Phi$ is not proper (not even if $Q$ is compact) and
Theorem \ref{theorem:convex} does not apply.  But the homogeneity of
$\Phi$ also implies that $\Delta(M)$ is equal to the cone on
$\Delta(U)$ for any neighbourhood $U$ of the zero section.  In view of
Theorem \ref{theorem:localcone} this means that
$\Delta(M)=\Delta_{(q,0)}$, where $q$ is any point in $Q$.
Furthermore, $\Phi(M)=-\Phi(M)$, so $\Delta(M)=\Delta(M)^*$.  The
symplectic slice $W$ at $(q,0)$ to the $K$-action on $M$ is equal to
$T^*V$, where $V$ is the slice $T_qQ/T_q(Kq)$ at $q$ to the $K$-action
on $Q$.  This implies that $W=V+JV=V\co$ for a suitable
$K_q$-invariant complex structure $J$ on $W$, and so the variety
$G\times^{(K_q)\co}W$ is the complexification of the real-algebraic
variety $K\times^{K_q}V$.  We have proved:

\begin{theorem}\label{theorem:cotangent}
For every connected $K$-manifold $Q$\upn, the set $\Delta(T^*Q)$ is a
rational convex polyhedral cone.   It is invariant under the involution
$*$ and equal to the momentum cone of the complexification of the
$K$-variety $K\times^{K_q}V$.  Here $V=T_qQ/T_q(Kq)$ is the slice at an
arbitrary point $q\in Q$.   
\qed
\end{theorem}

For instance, let $L$ be a closed subgroup of $K$ and let $Q$ be the
homogeneous space $K/L$.  Then the theorem says that
$\Delta\bigl(T^*(K/L)\bigr)=\Delta(G/L\co)$.

\subsection{Symplectic quotients}\label{subsection:reduction}

Let $M$ be a Hamiltonian $K$-manifold with momentum map $\Phi\colon
M\to\lie k^*$.  Let $L$ be a closed normal subgroup of $K$.  Then $M$
is a Hamiltonian $L$-space with momentum map
$\Phi_{(L)}=\iota^*\circ\Phi\colon M\to\lie l^*$, where $\iota$ is the
inclusion of $\lie l$ into $\lie k$.  Let $\bar K$ be the Lie group
$K/L$.  The kernel of $\iota^*$ can be identified in a natural way
with $\bar{\lie k}^*$, the dual of the Lie algebra of $\bar K$.  Let
$\mu$ be any point in $\lie z(\lie k)^*$, where $\lie z(\lie k)$
denotes the centre of $\lie k$.  The symplectic quotient of $M$ at the
level $\iota^*\mu\in\lie z(\lie l)^*$ with respect to the $L$-action,
$$
M_{\iota^*\mu}=M_{\iota^*\mu,L}=\Phi_{(L)}\inv(\iota^*\mu)/L
=\Phi\inv\bigl(\mu+\bar{\lie k}^*\bigr)\bigm/L,
$$ 
is a stratified Hamiltonian $\bar K$-space.  (Cf.\ \cite{sj:st}.)  A
momentum map $\Phi_{(\bar K)}\colon M_{\iota^*\mu}\to\bar{\lie
k}^*\subset\lie k^*$ for the $\bar K$-action on $M_{\iota^*\mu}$ is
induced by the map $\Phi\inv\bigl(\mu+\bar{\lie k}^*\bigr)\to\bar{\lie
k}^*$ sending $m$ to $\Phi(m)-\mu$.  (Up to a shift by an element of
$\lie z(\bar{\lie k})^*$, the map $\Phi_{(\bar K)}$ only depends on
the point $\iota^*\mu$.)  It is easy to calculate
$\Delta(M_{\iota^*\mu})$ in terms of $\Delta(M)$.  Let $\bar T$ be the
maximal torus $T/(T\cap L)$ of $\bar K$.  Then $\bar{\lie t}^*$ is
naturally isomorphic to $\lie t^*\cap\bar{\lie k}^*$, and the
intersection $\bar{\lie t}^*_+=\tplus\cap\bar{\lie k}^*$ is a Weyl
chamber of the pair $(\bar K,\bar T)$.  The following result is now
obvious (regardless of whether $M_{\iota^*\mu}$ is smooth or not).

\begin{proposition}\label{proposition:reduction}
$\Delta(M_{\iota^*\mu})=\bigl(-\mu+\Delta(M)\bigr)\cap\bar{\lie
k}^*$.  Therefore\upn, if $\Delta(M)$ is a closed convex polyhedral
subset of $\tplus$\upn, then $\Delta(M_{\iota^*\mu})$ is a closed
convex polyhedral subset of $\bar{\lie t}^*_+$.
\end{proposition}

For example, consider the Hamiltonian $K\times K$-space $T^*K\cong
K\times\lie k^*$ with momentum map $(k,\nu)\mapsto(k\nu,-\nu)$.  Let
$L$ be any closed subgroup of $K$.  The momentum map for the
restriction of the action to $K\times L$ is
$(k,\nu)\mapsto(k\nu,-\nu|_{\lie l})$.  The symplectic quotient of
$T^*K$ at level $0$ with respect to the normal subgroup $\{1\}\times
L\subset K\times L$ is isomorphic as a Hamiltonian $K$-space to the
cotangent bundle of the homogeneous space $K/L$.  Proposition
\ref{proposition:reduction} tells us that the momentum map image of
$T^*(K/L)$ is the set $K\{\,\nu:\nu|_{\lie l}=0\,\}$, and hence
\begin{equation}\label{equation:reduction}
\Delta\bigl(T^*(K/L)\bigr)=K\lie l^\circ\cap\tplus,
\end{equation}
where $\lie l^\circ$ is the annihilator of $\lie l$ in $\lie
k^*$.  (This can also be seen from the equality
$\Delta\bigl(T^*(K/L)\bigr) =\Delta(G/L\co)$ proven in Section
\ref{subsection:cotangent} and the Peter-Weyl Theorem.)  Consequently,
the set $K\lie l^\circ\cap\tplus$ is a rational convex polyhedral
cone.

Now assume that $K$ is semisimple and take $L$ to be the maximal torus
$T$.  Kostant's convexity theorem \cite{ko:on} implies that $K\lie
t^\circ\cap\tplus=\tplus$.  (Consider the natural projection
$\iota^*\colon\lie k^*\to\lie t^*$.  Take any $\mu\in\tplus$.  By
Kostant's theorem the set $\iota^*(K\mu)\subset\lie t^*$ is equal to
the convex hull of the Weyl group orbit through $\mu$, which contains
the origin in $\lie t^*$.  Therefore $K\mu\cap\lie
t^\circ=K\mu\cap\ker\iota^*$ is not empty.)  We conclude from
\eqref{equation:reduction} that $\Delta\bigl(T^*(K/T)\bigr)=\tplus$.

More generally, take $L$ to be the centralizer $K_\sigma$ of a wall
$\sigma$ of the Weyl chamber $\tplus$.  It is not difficult to see from
the root space decomposition of the pair $(\lie k,\lie t)$ that
$K\lie k_\sigma^\circ\cap\tplus$ contains the ray through every
dominant root that is not perpendicular to the wall $\sigma$.  This is
insufficient information to determine
$\Delta\bigl(T^*(K/K_\sigma)\bigr)$ in general, but if $K$ is e.g.\ of
type $\operatorname{B}_2$ or $\operatorname{G}_2$, then this implies
that $\Delta\bigl(T^*(K/K_\sigma)\bigr)=\tplus$ for any wall
$\sigma\neq\{0\}$.  If $K=\SU(n)$ and $\sigma$ is the one-dimensional
wall spanned by either $\pi_1$ or $\pi_{n-1}=\pi_1^*$, one can easily
calculate by hand that $\Delta\bigl(T^*(K/K_\sigma)\bigr)=K\lie
k_\sigma^\circ\cap\tplus$ is equal to the ray spanned by the maximal
root $\alpha_1+\cdots+\alpha_{n-1}=\pi_1+\pi_{n-1}$.

As another application of Proposition \ref{proposition:reduction}, I
now give a proof of Proposition \ref{proposition:vertex}.  For any
$\mu\in\tplus$ and for any closed subgroup $L$ of $K_\mu$, let $M$ be
the space $X(\mu,L,\{0\})$, the local model of Definition
\ref{definition:model} with trivial symplectic slice $W$.  By
\eqref{equation:bundle}, $M$ is the bundle over the coadjoint orbit
$K\mu$ with fibre $T^*(K_\mu/L)$ furnished with the minimal-coupling
form.  Put $\lie m=\lie k_\mu/\lie l$.  Then $\lie m^*$ is canonically
isomorphic to the annihilator of $\lie l$ inside $\lie k_\mu^*$.
Clearly, the point $m=(1,0,0)$ in $K\times^L\lie m^*\cong M$ has the
property that $K_m=L$ and $\Phi(m)=\mu$.  Also, the symplectic
cross-section of $M$ at $m$ is just the cotangent bundle
$T^*(K_\mu/L)$ with its standard momentum map shifted by $\mu\in\lie
z_\mu^*$.  By \eqref{equation:reduction}, the local momentum cone of
$M$ at $m$ is therefore equal to
\begin{equation}\label{equation:modelcone}
\Delta_m=\mu+\bigl(K_\mu\lie m^*\cap\lie t^*_{+,\mu}\bigr),
\end{equation}
where $\lie t^*_{+,\mu}$ denotes the positive Weyl chamber of $\lie
k_\mu$.  Now assume that $K_\mu=[K_\mu,K_\mu]L$, or, in other words,
$\lie k_\mu=[\lie k_\mu,\lie k_\mu]+\lie l$.  Because of the
decomposition $\lie k_\mu^*=[\lie k_\mu,\lie k_\mu]^*\oplus\lie
z_\mu^*$, this is equivalent to $\lie m^*\cap\lie z_\mu^*=\{0\}$, or:
\begin{equation}\label{equation:emptyintersection}
K_\mu\lie m^*\cap\lie z_\mu^*=\{0\}.
\end{equation}
Now the Weyl chamber $\lie t^*_{+,\mu}$ is the product of $\lie
z_\mu^*$ and the Weyl chamber of the semisimple part, $\lie
t^*_{+,\mu}\cap[\lie k_\mu,\lie k_\mu]$, which is a proper cone.  So
the cone $K_\mu\lie m^*\cap\lie t^*_{+,\mu}$ could only fail to be
proper if $K_\mu\lie m^*$ contained a nontrivial linear subspace of
$\lie z_\mu^*$.  But this is impossible because of
\eqref{equation:emptyintersection}.  By \eqref{equation:modelcone},
the point $\Phi(m)$ is therefore a vertex of $\Delta(M)$.  This
completes the proof of Proposition \ref{proposition:vertex}.


\providecommand{\bysame}{\leavevmode\hbox to3em{\hrulefill}\thinspace}


\end{document}